\documentclass[12pt]{article}     

\input epsf	

\date  {{\small Feb. 9, 1998; revised Sept. 14, 1998. }} 

\def\abst#1{\begin{minipage}{5.25in}
{\noindent   \normalsize
{\bf Abstract} #1}  \\ \end{minipage}  }

\def\be{\begin{equation}}
\def\ee{\end{equation}}
\def\bea{\begin{eqnarray}}
\def\eea{\end{eqnarray}}
\newcommand{\eq}[1]{eq.~(\ref{#1})} 

\def\C{{\mathcal C}}   
\def\CAP{\mathrm {Cap}}   
\def\E{{\mathrm E}}  
\def\F{{\mathcal F}}   
\def\P{Prob}   

\def\too#1{\parbox[t]{.4in} {$\longrightarrow\\[-9pt] {\scriptstyle #1}$}}
\def\dim{\overline{dim}_{B}}

\def\abs#1{\left\vert #1 \right\vert}
\def\blackbox{{\vrule height 1.3ex width 1.0ex depth -.2ex}\hskip 1.5
truecm}

\def\eps {\varepsilon}

\def\R  {{\bf R}}
\def\Z  {{\bf Z}}

\newcounter{masectionnumber}
\newcommand{\masect}[1]{\setcounter{equation}{0}
  \refstepcounter{masectionnumber} \vspace{1truecm plus 1cm} \noindent 
    {\large\bf \arabic{masectionnumber}. #1}\par \vspace{.2cm}
      \addcontentsline{toc}{section}{\arabic{masectionnumber}. #1}  
    }
 \renewcommand{\theequation}
    {\mbox{\arabic{masectionnumber}.\arabic{equation}}}
      
    \newcounter{masubsectionnumber}[masectionnumber]
\newcommand{\masubsect}[1]{
    \refstepcounter{masubsectionnumber} \vspace{.5cm} \noindent
  {\large\em \arabic{masectionnumber}.\alph{masubsectionnumber} #1}
    \par\vspace*{.2truecm}
 
\addcontentsline{toc}{subsection}          
 {\arabic{masectionnumber}.\alph{masubsectionnumber}\hspace{.1cm} #1} 
    }

\newtheorem{lem}{Lemma}[masectionnumber]
\newtheorem{thm}[lem]{Theorem}

\newcommand{\startappendix}{ \setcounter{masectionnumber}{0} 
 \renewcommand{\theequation}
    {\mbox{\Alph{masectionnumber}.\arabic{equation}}}          
 \renewcommand{\thelem}
    {\mbox{\Alph{masectionnumber}.\arabic{lem}}}          
\addcontentsline{toc}{section}{Appendix }
  } 
\newcommand{\maappendix}[1]{          
    \setcounter{equation}{0}
  \refstepcounter{masectionnumber} \vspace{1truecm plus 1cm} \noindent
    {\large\bf \Alph{masectionnumber}. #1}\par \vspace{.2cm} 
    
      \addcontentsline{toc}{section}{\Alph{masectionnumber}. #1 }   }
      
    \newcounter{masubapp}[masectionnumber]
\newcommand{\masubappendix}[1]{        
    \refstepcounter{masubapp} \vspace{.5cm} \noindent
  {\large\em \Alph{masectionnumber}.\alph{masubapp} #1}
    \par\vspace*{.2truecm}
 
\addcontentsline{toc}{subsection}          
 {\Alph{masectionnumber}.\alph{masubapp}\hspace{.1cm} #1 }     
      }

\begin{document}
\title{\vspace*{-.35in}
H\"older Regularity and Dimension Bounds for Random Curves}
\author{M. Aizenman 
${}^{(a)}$
\qquad and \qquad A. Burchard ${}^{(b)}$\\ \hskip 1cm
\vspace*{-0.05truein} \\
\normalsize \it  ${}^{(a)}$ School of Mathematics, 
Institute for Advanced Study, Princeton NJ 08540.\thanks{Permanent 
address:  Departments of Physics and  Mathematics, Jadwin Hall, 
Princeton University, P. O. Box 708, Princeton, NJ 08540.} \\ 
\normalsize \it ${}^{(b)}$ Department of Mathematics, 
Princeton University,  Princeton NJ 08544.}
\maketitle
\thispagestyle{empty}        

   \begin{minipage}{5.25in}
     { \em 
     \begin{center}
     Dedicated to the memory of \/  {\em Roland L. Dobrushin}:  \\
     Exemplary in science and in life.    
      \end{center}    }
    \end{minipage}   

\bigskip \bigskip

\abst{Random systems of curves exhibiting fluctuating features on 
arbitrarily small scales ($\delta$) are often encountered in 
critical models.  For such systems it is shown that 
scale-invariant bounds on the probabilities of crossing 
events imply that typically all the realized curves 
admit H\"older continuous parametrizations with a common 
exponent and a common random pre\-factor, which in the scaling 
limit ($\delta\to 0$) remains stochastically bounded. 
The regularity is used for the construction of scaling limits, 
formulated in terms of probability 
measures on the space of closed sets of curves.  
Under the hypotheses presented here the limiting measures are 
supported on sets of curves which are H\"older continuous 
but not rectifiable, and have Hausdorff dimensions strictly 
greater than one.  The hypotheses are known to 
be satisfied in certain two dimensional percolation models. 
Other potential applications are also mentioned.    }

\vfill

\bigskip  \bigskip \bigskip
 
\newpage
\begin{minipage}[t]{\textwidth}
\tableofcontents
	\end{minipage}
	
\newpage
\masect{Introduction}  
\label{sect:intro}
\vskip-0.6cm

\masubsect{General framework}

We consider here curves in $R^{d}$ which are shaped  on
many scales, in a manner found in various critical models. 
The framework for the discussion are random systems, where the
random object is expressed as a closed  collection of polygonal curves
of a small step size $\delta$.  Our main results are
general criteria for establishing that, as the short distance cutoff is  
sent to zero,  the curves retain
a certain degree  of regularity, yet at 
the same time are intrinsically rough.
In some cases the object of study is a single random curve. In others,
the random object contains many curves; in such situations, 
the regularity estimates are intended to cover the entire collection.

The criteria developed here can be applied to various stochastic 
geometric models.  
In the Appendix we mention as examples critical percolation, 
the minimal spanning trees in random geometry, 
the frontier of two-dimensional Brownian motion,  
and the level sets of a two-dimensional random field.  

While our discussion does not require familiarity with any of 
these examples, let us comment that an important feature 
they share is the existence of two very different length scales: 
the {\em microscopic} scale on which  the model's building variables 
reside, and the {\em macroscopic} scale on which the connected curves 
are tracked.   In such situations it is natural to seek a meaningful 
formulation for the {\em scaling limit}, at which the 
microscopic scale ($\delta$) is taken to zero.  
The regularity established here enables such a construction 
through compactness  arguments.

To introduce the results let us start with some of the terminology.
\begin{itemize}
\item[i.]
We denote by ${\cal S}_\Lambda$  the space of curves in a closed subset 
$\Lambda\subset \R^d$, 
with the metric defined in Section~\ref{sect:curves}. \\
The symbol $\C$ is reserved here for individual curves, and 
$\F$ for  sets of curves.    
The space of closed sets of curves in $\Lambda$ is 
denoted $\Omega_\Lambda$.  
\item[ii.]  
A {\em configuration of curves}  in $\R^d$ 
with a short-\-distance cutoff $\delta \in (0,1]$  is a 
collection of polygonal paths of step size $\delta$ which forms 
a closed subset $\F_{\delta}(\subset {\cal S}_\Lambda)$.   
\item[iii.]
A {\em system of random curves} with varying short-\-distance 
cutoff is described by a collection of probability measures 
$\{ \mu_{\delta}(d\F_{\delta}) \}_{0< \delta \le \delta_{max} } $  
on $\Omega_\Lambda$, where each $\mu_{\delta}$ describes random 
sets  of curves in $\Lambda$ consisting of polygonal paths of 
step size $\delta$.  (We shall often take $\delta_{max} = 1$).  
\end{itemize}
 
To summarize: the individual realizations of the random 
systems are closed sets of curves denoted $\F_\delta(\omega)$.  
The entire system will occasionally 
also be referred to by the symbol~$\F$.

\begin{figure}[t]
      \begin{center} 
      \leavevmode
  	\epsfysize=2.5in
      \epsfbox{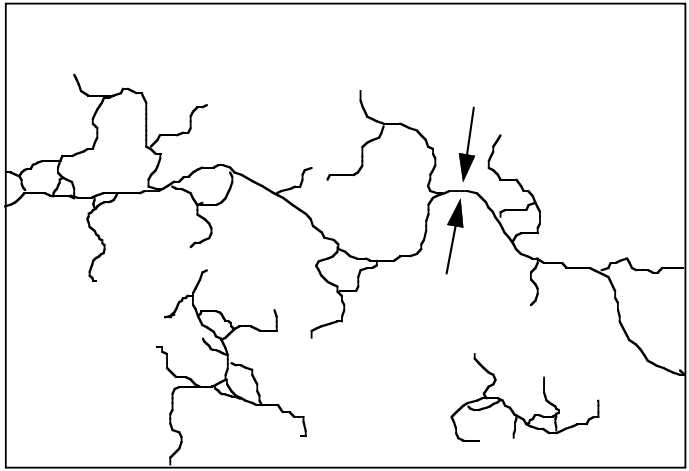}
  \caption{\footnotesize  Schematic depiction of some of the realized 
  paths in a critical percolation model.  The model is formulated 
  (possibly on a lattice)  on a scale which is much 
  smaller than the length of the depicted region.  
   Typical configurations exhibit large-scale connected clusters, but 
   the connections are tenuous.  In this work we discuss the 
   regularity properties of the self-avoiding paths 
   supported on such clusters. }
  \label{fig:chokepoint}
  \end{center}
  \end{figure}
  
We are particularly interested in 
statements concerning the probability measures 
$\{ \mu_{\delta}(d\F_{\delta})\}_{0< \delta \le \delta_{max}}$ 
which hold uniformly  in $\delta$, 
and thus provide information about the scaling limit. 
The following notion facilitates the formulation of such  
statements.   

\bigskip\noindent{\bf Definition} \ 
A random variable $X$ associated with a system $\F$ is said to 
be {\em stochastically bounded}, for $\delta \to 0$, if:  
\begin{itemize} \item[i.] a version $X_\delta$ is defined for each 
$0< \delta (\le \delta_{max})$, and 
\item[ii.] for every $\eps > 0$ there is $u  < \infty$ 
such that 
\be 
Prob_{\delta}\left( |X_{\delta}(\omega)| \ge u \right) \ \le \ \eps 
\ee 
uniformly in $\delta$. 
\end{itemize}
A random variable is said to be {\em stochastically bounded away from 
zero} if its inverse is stochastically bounded.   
\bigskip
 
%

\masubsect{Main results}
 
The main theorems in this paper are derived under two hypotheses, which 
require 
scale-invariant bounds on certain crossing events.  
In order to apply these results, the hypotheses
need to be verified using specific information on the given system. 
The conditions are known to be satisfied in certain   two dimensional critical 
percolation models, through the ``Russo-Seymour-Welsh theory'' 
\cite{R,SW} and its recent extensions \cite{Alex-RSW}.  They are 
expected to be true for critical percolation also in higher dimensions, 
though not for $d>6$~\cite{Aiz-ISC}.  Another proven example is the 
minimal spanning tree in two dimensions~\cite{ABNW}, where the 
criticality is {\em self-induced}  in the sense of~\cite{BTW}.  
These, and other systems, are presented in the Appendix.  

 \begin{figure}[t]
      \begin{center} 
      \leavevmode
  	\epsfysize=2in
      \epsfbox{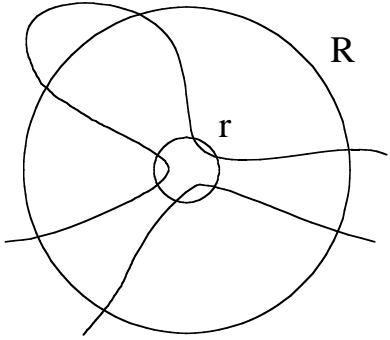}
  \caption{\footnotesize  
   A 6-fold crossing of an annulus $D(x;r,R)$.  
  Under the hypothesis {\bf H1} the probability of a
  k-fold crossing of an annulus $D(x;r,R)$ is
  bounded by $K_k (r/R)^{\lambda(k)}$.  The regularity estimates are 
  based on the absence of k-crossings for $k$ large enough and  
  $0 < r <<  R= r^{1-\eps} << 1$. }
  \label{fig:6crossing}
  \end{center}
  \end{figure}

The first condition concerns repeated crossings of the 
spherical shells 
\be
D(x; r, R) = 
\left\{y \in \R^{d} \ {\big |}\ r \le |y-x| \le R \right\} \; .
\label{eq:def-shell}
\ee
The assumption is:
\begin{itemize}
\item[{\bf (H1)}] {\em Power-bound on the probability of  multiple
crossings.}  For all $k < \infty$ and  for all spherical shells 
with radii $\delta\le r\le R\le 1$, the following bound holds 
uniformly \\   
in $\delta$ 
\be 
Prob_{\delta}\left( \begin{array}{c} 
       \mbox{$D(x;r,R)$  is traversed by $k$ separate} \\
\mbox{ segments of a curve in $\F_{\delta}(\omega)$}
\end{array} \right) \ \le\ K_k \left(\frac{r}{R}\right)^{\lambda(k)}
\label{eq:suff-upper-1}
\ee
with  some $K_k\le\infty$ and 
\be 
 \lambda(k)\ \to \infty \ \qquad  \mbox{as $k \to \infty$} \;  .
\label{eq:suff-upper-2}
\ee
\end{itemize}

Based on this assumption, we derive the following result.  

\begin{thm} {\rm (Regularity)}
Let $\F$ be a system of random curves in a compact region 
$\Lambda \subset \R^d$, with variable short-distance cutoff 
$\delta>0$, and assume the hypothesis $\bf H1$ is satisfied.  
Then for each $\eps > 0$ 
all the curves $\C \in \F_{\delta}(\omega)$ 
can be simultaneously parametrized by continuous functions
 $\gamma(t)$ $t\in [0,1]$, 
such that for  all  $0\le t_{1} < t_{2} \le 1$:
\be
| \gamma (t_{1}) - \gamma (t_{2}) | \ 
 \le \kappa_{\eps;\delta}(\omega)\ g(\mbox{diam\/}(\C) )^{1+\eps} \  
  \ |t_{1}-t_{2}|^{1 \over d - \lambda(1) + \eps}\; ,
\label{eq:holder-gen}
\ee
with a random variable $\kappa_{\eps;\delta}(\omega)$ 
(common to all $\C\in \F_{\delta}(\omega)$) which stays 
stochastically bounded as $\delta \to 0$.  The second factor 
depends on the curve's diameter through the function 
\be
g(r) \ = \ r^{-{\lambda(1) \over d - \lambda(1)} }\;  . 
\ee 
\label{thm:regularity}
\end{thm} 

\noindent {\bf Remarks:} \
\noindent {\em i.\/} 
The conclusion of Theorem~\ref{thm:regularity} is stated 
in terms of the existence of H\"older continuous parametrizations, 
which offers a familiar criterion for regularity.
We actually find it convenient to develop the results in terms of 
{\em tortuosity bounds}, that is, upper bounds on the function 
$M(\C,\ell)$ defined as the minimal number of 
segments produced if the curve $\C$ is partitioned into 
segments of diameter no greater that $\ell$. 
The two notions are linked in Section~\ref{sect:curves}.     

{\em ii.\/} The dependence of the H\"older constant
in \eq{eq:holder-gen} on the diameter of
the curve  can be removed by lowering the
H\" older exponent below $1/d$; indeed, interpolating 
between~\eq{eq:holder-gen} and the trivial relation
$| \gamma (t_{1}) - \gamma (t_{2})| \ \le\ \mbox{diam}(\C)$
gives
\be
| \gamma (t_{1}) - \gamma (t_{2}) | \ 
 \le \left(\kappa_{\eps;\delta}(\omega)\right)^{1+\tilde\eps-d/\lambda}
  \ |t_{1}-t_{2}|^{1 \over d + \tilde\eps}\; ,
\label{eq:holder-unif} 
\ee
where $\tilde \eps$ is small when $\eps$ is small.

{\em iii.\/}  For the main conclusion that the  
curves retain H\"older continuity at {\em some} $\alpha > 0$, 
it suffices to require instead  of {\bf H1} that 
\eq{eq:suff-upper-1} holds for some $k$ with  
$\lambda(k) > d$.   While this is clearly a weaker assumption 
than  \eq{eq:suff-upper-2}, so far it has been proven only in 
situations in which also (\ref{eq:suff-upper-2}) holds.  

In addition to being of intrinsic interest, 
the above regularity property permits to construct the scaling 
(continuum) limit.  The basic question concerning this limit is:

{\em  \begin{enumerate}
\item[Q1.] Is the collection of probability measures 
$\{\mu_{\delta}(d \F_{\delta})\}_{\delta}$ tight?  
\end{enumerate}}
Tightness means that, 
up to remainders that can be made arbitrarily small, 
the measures $\mu_{\delta}$ 
share a common compact support.  A positive answer 
to {\em Q1} implies the 
existence of limits $ \lim \mu_{\delta, p(\delta)}$ at least along
some sequences of $\delta_n \to 0$ (\cite{Bill}).  
Without tightness, one cannot rule out 
the possibility that, as the cutoff is removed, curves give way
to more general continua.  The range of possibilities 
is rather vast: curves can converge (in a weaker sense than used
here) to continua which do not support any continuous curve.  

Theorem~\ref{thm:regularity} yields a positive answer to 
{\em Q1}, as is explained here in  Section~\ref{sect:compact}. 
There are a number of dimension-like quantifiers for the 
description of the  curves emerging in the scaling limit.
Among them are (see Section~\ref{sect:curves}) 
\begin{itemize}
\item
the Hausdorff dimension  $dim_{\cal H}\C$,
\item
the upper box dimension (also known as the Minkowski dimension)
$\dim{\C}$,
\item
and the reciprocal of the optimal H\"older regularity exponent 
\be
\alpha(\C) \ = \ \sup\left\{ \alpha  | 
\begin{array}{c}
 \mbox{ $\C$ can be parametrized as $\{ \gamma(t)\}_{0\le t\le 
 1}$  with}  \\
 \mbox{$\abs{\gamma(t) - \gamma(t')}  \ \le \ K_{\alpha}
 \abs{t - t'}^{\alpha}\; $,  for all  $0 \le t \le t' \le 1$}	
\end{array}
\right\}  \; .
\ee
\end{itemize}

The following result is derived in Section~\ref{sect:compact}.

\begin{thm} {\rm (Scaling limit)}\ 
For any system $\F$ of random curves in a compact set 
$\Lambda\subset \R^d$, hypothesis {\bf H1} implies  that
the limit
\be 
\lim_{n\to \infty} \mu_{\delta_n}(d\F) \ =:\  \mu(d\F)
\label{eq:limit}
\ee 
exists at least for some  sequence $\delta_{n} \to 0$.
The limiting probability measure (on $\Omega_\Lambda$) 
is supported on configurations $\F$ containing only paths 
with 
\be 
\dim{\C} \ = \ \alpha(\C)^{-1} \ \le \ d - \lambda(1) 
\label{eq:minkowski}
\ee
and
\be
\mbox{dim\/}_{\cal H}(\C) \ \le d - \lambda(2) \; .
\label{eq:bbdim}
\ee 
\label{thm:scaling}
\end{thm}
\noindent (The improvement in the dimension
estimate of the last inequality over the preceding one is  based on 
considerations of the {\em ``backbone''}.)

The sense of convergence in \eq{eq:limit} 
can be expressed by 
saying that there exists a family of {\em couplings} 
consisting of probability measures $\rho_n(d\F_{\delta_n}, d\F)$ 
such that:
\begin{itemize}
\item[i)]
 the marginal distributions satisfy
\be
 \rho_n(d\F_{\delta_n}, \Omega_{\Lambda}) \ = \  
 \mu_{\delta_n}(d\F_{\delta_n}) 
 \; ,  \qquad \;   \rho_n(\Omega_{\Lambda}, d\F) \ = \  \mu(d\F) 
\ee 
\end{itemize}
and 
\begin{itemize}
\item[ii)]    the two components are close, in the sense that 
\be
\int_{ \Omega_{\Lambda}\times  \Omega_{\Lambda}}  dist(\F_{\delta_n},\F) 
\; \;  \rho_n(d\F_{\delta_n}, d\F)  \  \too{n\to \infty} \ 0 \; ,
\label{eq:vasershtein}
\ee 
\end{itemize} 
with the distance between two configurations of curves defined by
the Hausdorff metric: 
\be 
dist(\F,\F') \le \eps \quad  \Leftrightarrow \quad 
\left\{ 
\begin{array}{l}
	\mbox{for every $\gamma \in \F$ there is  $\gamma' \in \F'$}  \\
	\mbox{with $\sup_{t}|\gamma(t) - \gamma'(t)| \le \eps$} \\ 
	\mbox{ \quad and {\em vice-versa} ($\F\leftrightarrow \F'$)  }  
	\end{array}  \right.   \; . 
\ee 

The positive answer to {\em Q1}  invites a number of other 
questions, such as:
{\em 
\begin{enumerate} 
\item[Q2.] Is the limit independent of the sequence 
$\{\delta_n\}$, 
and is it shared by other models 
with different short-scale structure?
\end{enumerate}
}
This question is beyond the scope of the present work.  
In some of the models of interest it is related to the purported 
universality of critical behavior.

\begin{figure}[t]
      \begin{center} 
      \leavevmode
  	\epsfysize=2in
      \epsfbox{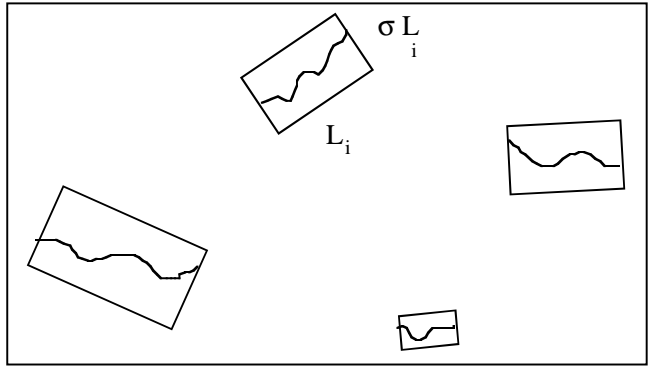}
  \caption{\footnotesize  A simultaneous crossing event for a family 
  of cylinders of common aspect ratio.  In hypothesis {\bf H2} 
  the probability of such an events is assumed to be less than 
  $Const. \rho^k$, with $\rho < 1$  (here $k=4$).  
  The implication is a uniform lower bound 
  on the Hausdorff dimensions of all the curves in the configuration.}
  \label{fig:H2}
  \end{center}
  \end{figure}
  
To establish some minimal roughness for the realized curves, 
we require a hypothesis on the probability of 
simultaneous crossings of a family of cylinders.  
It suffices to restrict the assumption to spatially
separated cylinders, with much latitude in the exact 
definition of the term.   

\bigskip\noindent{\bf Definition} \ 
A collection of sets $\{A_j\}$ is {\em well separated}
if the distance of each set $ A_j $ to 
the other sets $\{A_i\}_{i\not = j}$
is at least as large as the diameter of $A_j$.   
\bigskip

The relevant hypothesis is:
\begin{itemize}
\item[{\bf (H2)}] {\em Power bound on the probability of simultaneous 
crossings.} \ There exist a cross section $\sigma> 0$ 
and some $\rho < 1$ with which for every collection of $k$ well 
separated cylinders, $A_1,\ldots, A_k $, of aspect ratio $\sigma$ 
and lengths $\ell_1,\ldots, \ell_k \ge \delta$ 
\be
\hbox{Prob}_\delta\left( \begin{array}{c}
  \mbox{each $A_j$ is traversed (in the} \\ \mbox{ long 
  direction) by a curve in $\F_{\delta}$   }
  \end{array} 
\right)\ \le\ K \rho^k \; . 
\label{eq:suff-lower}
\ee
\end{itemize}
\bigskip

An effective way to express curve roughness is in terms of   
{\em capacity} lower bounds.  The capacity ${\CAP}_{s;\ell}(A)$ 
of a set $A\subset \R^d$ is defined in Section~\ref{sect:sparse}.  
For the purpose of this summary it suffices to note that the capacity 
of a fixed set $A$ increases with the parameter $\ell$, 
and that it  provides the following  lower bounds on coverings. 
\begin{itemize}
\item[i.]  
For every covering of $A$ by sets $\{B_j \}$ of diameter at least 
$\ell$,
\be 
\sum_{j} \left( \mbox{diam} B_j \right)^s \ \ge \ 
{\CAP}_{s;\delta}(A) \ .
\label{eq:cap-haus1}
\ee 
\item [ii.]
The minimal number of elements for a covering of $A$ by 
sets of diameter $\ell$ satisfies
\be
N(\C,\ell)\ \ge\     \CAP_{s;\ell}(A) \cdot \ell^{-s}
      \; .
\label{eq:cap-N}
\ee
\item[iii.]   
The behavior of the capacity  for small $\ell$ 
provides information on the {\em Hausdorff dimension}:
\be 
\inf_{0< \ell \le 1 } {\CAP}_{s;\ell}(A)  \ > \ 0 \ 
 \quad  \ \Rightarrow \quad dim_{\cal H}(A) \ge s \; .
\ee 
\end{itemize}
(Proof of {\em i.\/} 
is given in Section~\ref{sect:sparse}; {\em ii.\/} and {\em iii.\/} 
are direct consequences.)

\begin{thm} {\rm (Roughness)  }
If a system $\F$ of random curves in a compact subset $\Lambda\subset 
\R^d$ 
with a variable short distance cutoff 
satisfies the hypothesis {\bf H2}, then 
there exists $d_{min} > 1$ such that for any fixed $r > 0$  and $ s > 
d_{min}$
the random variable
\be 
T_{s,r;\delta}(\omega) \ := \ \inf_{
 \C \in \F: \ diam({\cal C}) \ge r }    {\CAP}_{s;\delta}(\C) 
\label{eq:roughness}
\ee
stays stochastically bounded away from zero, as $\delta \to 0$.  

Furthermore, any scaling limit  of the measures $\mu_{\delta}$,  
$\mu=\lim_{\delta_{n}\to 0}\mu_{\delta_n}$,
is supported on  configurations containing only curves with
\be
dim_{\cal H}\C \ \ge \ d_{min} \ (>1) \; .
\ee 
In particular, the scaling 
limit contains no rectifiable curves.
\label{thm:roughness}
\end{thm}

Let us  note that Theorem~\ref{thm:roughness} complements  
Theorem~\ref{thm:regularity}, since by the monotonicity properties
of the capacity, we can combine 
\eq{eq:roughness} with \eq{eq:cap-N}  to obtain
\be
N(\C, \ell) \ \ge \  {\CAP}_{s,\ell}(\C)\cdot \ell^{-s}\ \ge\ 
T_{s,r;\delta}(\omega) \ \ell^{-s}  \; 
\label{eq:min}
\ee 
for all $\ell\ge\delta$ 
whereas under the condition  (\ref{eq:holder-gen}) 
\be 
N(\C, \ell) \ \le \  
\left[{\kappa_{\eps;\delta}(\omega)\ 
g(\mbox{diam\/}\C)^{1+\eps} \over \ell }\right]^{d -\lambda(1) + \eps} \ 
 \; .
\label{eq:holder-N}
\ee 

In particular, \eq{eq:min}  implies 
that the minimal number of steps of the lattice size ($\delta$) 
needed in order 
to advance distance $L$ exceeds $Const. \ (L/ \delta)^\tau$.  
Some bounds of this form were previously obtained  
for the ``lowest path'' in 2D critical 
percolation models in the work of  Kesten and Zhang \cite{KZ93}, 
who refer to the optimal value of $\tau$ as the {\em tortuosity}  
exponent. 
We slightly modify their terminology, by requiring the power bounds
to hold simultaneously on all scales.    

The assumption in Theorem~\ref{thm:roughness} 
can be weakened by restricting  \eq{eq:suff-lower}  
to collections of cylinders of comparable dimensions, but then the 
conclusion will be stated in terms of the box dimension.
 
In models where spatially separated events 
are independent, condition {\bf H2} is implied by {\bf H1}, provided 
the parameter $\lambda(1)$ of \eq{eq:suff-upper-2} is positive.  
A similar observation applies to models without strict 
independence but with a correlation length of only microscopic size, 
such as the droplet percolation model.  

There is a considerable disparity between the upper and the lower 
bounds derived here for the dimensions of curves in the scaling limit. 
Part of the reason is that our lower bounds are far from sharp.  
However, we also expect some of the systems considered here 
(e.g. the percolation models) to exhibit 
simultaneously curves of different dimensions.  

The organization of the paper is as follows.   
In Section~\ref{sect:curves} we prepare for the discussion of 
random systems by clarifying some notions pertaining to single curves.
Introduced there is the concept of {\em tortuosity}, 
which provides a measure of roughness manifestly independent of 
parametrization. The associated tortuosity exponent coincides with 
Richardson's exponent $D$.  It is related here  to 
the degree of H\"older regularity achievable through 
reparametrization (Theorem~\ref{thm:alphatau}).  
Moreover, under  the {\em tempered crossing} condition 
the tortuosity exponent coincides with the curve's upper 
box dimension (Theorem~\ref{thm:taubox}).  
In Section~\ref{sect:random}  we apply these relations to general
systems of random curves, and prove the regularity result,
Theorem~\ref{thm:regularity}. To tighten the regularity estimate 
we briefly discuss the concept of the {\em backbone}.  
Section~\ref{sect:compact}  deals with the 
construction of scaling limits and the proof
of Theorem~\ref{thm:scaling}, based on 
the afore-mentioned regularity properties.  
The proof of the roughness lower bounds is split into two parts.
In Section~\ref{sect:sparse} we derive a deterministic statement 
(Theorem~\ref{thm:lower-det})
which presents a criterion for the roughness
of a curve,  based on the assumption that 
the {\em straight runs} of the curve are {\em sparse}.
The analysis exploits the relation of the dimension with {\em capacity}, 
and involves suitable energy estimates.  
In Section~\ref{sect:lower}, we apply
this result to random systems, and prove Theorem~\ref{thm:roughness} 
by establishing that {\bf H2} implies the sparsity of straight runs.   
The Appendix includes examples of systems for which the general  
Theorems yield results of interest within the specific context.   

%
%
\masect{Analysis of curves through tortuosity}
\label{sect:curves}

In this section we introduce the space of curves, and the 
notion of tortuosity.  
The two basic results are Theorem~\ref{thm:alphatau}, which 
relates the tortuosity exponent with 
the optimal H\"older continuity exponent, and 
Theorem~\ref{thm:taubox} which provides useful conditions under 
which the tortuosity exponent agrees with the upper box dimension, 
and thus is finite.

\masubsect{The space of curves}

We regard curves as equivalence classes of continuous 
functions, modulo re\-para\-me\-tri\-za\-tions. More precisely, 
two continuous functions $f_1$ and $f_2$ from
the unit interval into $\R^d$ describe the
same curve if and only if there exist two 
monotone continuous bijections $\phi_i : [0,1] \rightarrow [0,1] $, 
$i=1,2$,  so that 
$f_1\circ\phi_1=f_2\circ\phi_2$.  

Recall that the space of curves in 
a closed subset $\Lambda\subset \R^d$ is denoted here by $S_\Lambda$. 
The distance between two curves is measured by:
\be
{\rm d}(\C_1,\C_2)\ :=\ 
\inf_{\phi_1,\phi_2}\ \sup_ {t\in [0,1]}
\abs{f_1(\phi_1(t))-f_2(\phi_2(t))}\; ,
\label{eq:def-metric}
\ee
where  $f_1$ and $f_2$ is any pair of continuous functions representing
$\C_1$ and $\C_2$, and the infimum is over the set of all  strictly
monotone continuous functions from the unit interval  onto itself.

\begin{lem} Equation~(\ref{eq:def-metric}) defines a metric 
on the space of curves.
\label{lem:metric}
\end{lem}
\noindent{\bf Proof:} \ Clearly, ${\rm d}(\C_1,\C_2)$ is 
nonnegative, symmetric, satisfies the triangle inequality and 
${\rm d}(\C,\C)=0$.  To prove strict positivity, assume 
${\rm d}(\C_1,\C_2)=0$, and choose parametrizations $f_1$ and $f_2$.
We need to show that $f_1$ and $f_2$ describe the same curve, i.e.,
$\C_1=\C_2$. We may choose $f_1$ and $f_2$ to
be non-constant  on any interval.  Under these assumptions, there exist 
sequences of reparametrizations $\phi_1^i$ and $\phi_2^i$ such that 
\be
\sup_{t \in [0,1]} 
\abs{ f_1\circ \phi_1^i\circ (\phi_2^i)^{-1}(t) - f_2(t) }
\ =\ 
\sup_{t \in [0,1]} 
\abs{ f_1\circ \phi_1^i(t) - f_2\circ\phi_2^i(t) } 
\ \too{i\to\infty} 0\ ; .
\ee
Monotonicity and uniform boundedness imply (Helly's theorem)
that there are subsequences (again denoted $\phi_1^i$ and $\phi_2^i$)
so that $\phi_2^i\circ (\phi_1^i)^{-1}$ and their
inverses $\phi_1^i\circ (\phi_2^i)^{-1}$ converge
pointwise, at all but countably many points, to monotone 
limiting functions $\phi$ and $\tilde \phi$,
with   $f_1 = f_2\circ\phi$ and $f_2=f_1\circ\tilde\phi$.
To see that  $\phi$ has no discontinuities, note that 
jumps of $\phi$ would correspond to intervals where 
$\tilde \phi$ is is constant. But $\tilde \phi$ cannot be 
constant on an interval, since, by our choice of
parametrization, $f_2$ is not constant on any
interval.
\hfill\blackbox

\bigskip With this metric, ${\cal S}_{\Lambda}$ is complete but,
even for compact $\Lambda$ it is not 
compact.  This reflects the properties of 
the space of continuous functions  $C([0,1], \Lambda)$.

\masubsect{Measures of curve roughness}

Let $M(\C,\ell)$  be the  minimal number of segments needed for
a partition of a curve $\C$ into segments of diameter no greater 
than $\ell$.  Any bound on $M(\C,\ell)$ will be called a
{\em tortuosity bound}. 
In particular, we are interested in power bounds of the form
\be
M(\C,\ell) \ \le \ K_s \,\ell^{-s}  \; .
\label{eq:Ms}
\ee
Optimization over the exponents yields the following dimension-like 
quantity.

\bigskip\noindent{\bf Definition} \ {\it
For a given curve $\C$, 
\be
\tau(\C) \ = \ \inf\{ s > 0 | \  \ell^s \,M(\C,\ell) \ \too{\ell \to 0} 0
\} 
\label{eq:def-tau}
\ee
is called the {\em tortuosity exponent}.  
}\bigskip

There are a number of other ways of dividing a curve to short segments
which yield comparable results.  
Of particular interest to us is the observation that  
the tortuosity exponent can also be based on 
$\tilde M(\C,\ell)$, which we define
as the maximal number of points that can be 
placed on the curve so that successive points have distance
at least $\ell$. $M(\C,\ell)$ and $\tilde M(\C,\ell)$
are comparable, but have different continuity properties.

\begin{lem} $M(\C,\ell)$ and $\tilde M(\C,\ell)$ are
related by the inequalities
\be
M(\C,3\ell)\ \le\ \tilde M(\C,\ell)\ \le \ \inf_\eps M(\C,\ell-\eps)
\; .
\label{eq:claim-semicontinuous}
\ee
Furthermore, $M(\C,\ell)$ is lower semicontinuous, and $\tilde 
M(\C,\ell)$ is upper semicontinuous on the space of curves.
\label{lem:semicontinuous}
\end{lem}

\noindent{\bf Proof} \ The first inequality holds because
a segment of the curve of diameter  at least $3\ell$ certainly contains  
a point
that has a distance of at least $\ell$ from both
endpoints. The second inequality holds because no
segment of diameter less than $\ell$ can contain two points
of distance $\ell$ or more. The continuity properties 
follow easily from the fact that $M$ was defined through 
minimization and $\tilde M$ through maximization of cut points. 
\hfill\blackbox

It follows from Lemma~\ref{lem:semicontinuous} 
that the tortuosity exponent coincides 
with Richardson's exponent $D$ \cite{Richardson,Man}, 
which in ref.~\cite{DDC} was termed the ``divider dimension''.  
It was pointed out there that $D$ can take arbitrarily large 
values.  

From a different perspective, the curve's regularity may be 
expressed through the degree of {\em  H\"older continuity} 
achievable through
reparametrization.   One attempts to describe 
the curve
by means of a continuous function 
$\C = \{\gamma (t)\}_{0\le t \le 1}$,
satisfying:
\be
\abs{\gamma(t_1) - \gamma(t_2)}  \ \le \ K_{\alpha}
 \abs{t_1 - t_2}^{\alpha}
\; \mbox{ for all $0 \le t_1 \le t_2 \le 1$},
\label{def:holder}
\ee
with some exponent $\alpha >0$.
Greater values of the exponent correspond to higher degrees of
regularity, and thus one is interested in
\be
\alpha(\C) \ = \ \sup\left\{ \alpha |
\mbox{ $\C$ admits a parametrization satisfying \eq{def:holder}
with exponent $\alpha$}
\right\}  \; .
\ee

The tortuosity exponent may remind  one of the {\em upper box dimension}, 
which has a similar definition. Let $N(\C,\ell)$ be the minimal number 
of sets of
diameter $\ell$ needed to cover the curve. Then
\be
\dim({\C})\ :=\
\inf\{ s > 0 | \  \ell^s \,N(\C,\ell) \ \too{\ell \to 0} 0 \} \ .
\label{eq:def-box}
\ee
The two definitions are different, since a single set
of diameter $\ell$ may contain a large number of segments
of the curve. The box dimension can be calculated using only 
coverings with boxes taken from subdivisions of a fixed grid.

A trivial relation between the three parameters is
\be
\dim(\C)\ \le\ \tau(\C)\ \le\ \alpha(\C)^{-1}\; ,
\label{ineq:dim-trivial}
\ee 
which  follows immediately from
\be
N(\C,\ell)\ \le M(\C,\ell)\ \le\ 
\left\lceil\frac{K_\alpha}{\ell}\right\rceil^{1/\alpha}\ ,
\label{ineq:dim-trivial-singlescale}
\ee
where $\lceil x\rceil$ denotes the smallest integer at least as 
large as $x\,$.

\masubsect{Tortuosity and H\"older continuity}

It turns out  that the tortuosity exponent and the optimal 
H\"older exponent are directly related: 

\begin{thm} For any curve $\C$ in ${\cal S}_{\R^d}$,
\be 
\tau(\C)\ = \ \alpha(\C)^{-1}\; .  
\label{eq:alphatau}
\ee
\label{thm:alphatau}
\end{thm}

More explicitly, uniform continuity is equivalent to
a uniform upper tortuosity bound, as expressed in the following 
lemma.
 
\begin{lem} 
If a curve $\C$  in $\R^d$ admits a parametrization 
as  $\{\gamma(t)\}_{0\le t\le 1}$
so that for all $t_1$, t$_2$  in the unit interval
\be
\psi(\abs{\gamma(t_1)-\gamma(t_2)})\ \le\ \abs{t_1-t_2}\; ,
\label{ineq:Delta-gamma-1}
\ee
where $\psi:(0,1]\to(0,1]$ is a nondecreasing function,
then, for all $\ell\le 1$,
\be
\ M(\C,\ell) \ \le \ \left\lceil\frac{1}{\psi(\ell)}\right\rceil\; .
\label{ineq:M(C,l)-ceil}
\ee
Conversely, if 
\be
\ M(\C,\ell) \ \le \ \frac{1}{\psi(\ell)}
\label{ineq:M(C,l)}
\ee
for all $\ell\le 1$, then $\C$ can be parametrized as 
$\{\gamma(t)\}_{0\let\le1}$ with a function satisfying 
\be
\tilde\psi(\abs{\gamma(t_1)-\gamma(t_2)})\ \le\ 
\abs{t_1-t_2}\; ,
\label{ineq:Delta-gamma-2}
\ee
for all $0\le t_1 < t_2\le 1$, with 
\be
\tilde \psi(\ell)\ =\ 
{ \psi(\ell /2) \over 2(\log_2 (4/\ell))^2 } \; .
\label{eq:uniform-log} 
\ee
\label{lem:uniform}
\end{lem}

\noindent {\bf Proof:} The tortuosity bound in~\eq{ineq:M(C,l)-ceil}
follows from the uniform continuity condition 
in~\eq{ineq:Delta-gamma-1} with the definition of $M(\C,\ell)$,
by partitioning the curve into segments corresponding to time 
intervals of length $\psi(\ell)$.
To prove the reverse implication,  we need to construct a
parametrization satisfying the uniform continuity
condition in~\eq{ineq:Delta-gamma-2}, given that \eq{ineq:M(C,l)}
holds for $0<\ell\le 1$.
Choose an auxiliary parametrization of the curve, $\C=\{\gamma(s)\}$,
which is not constant on any interval. 
We associate with each curve
segment, $\C_s = \gamma([0,s])$, the {\em ``time of travel''}
\be
t_\eps(s)\ :=\ 
      \frac{\sum_n (n+1)^{-2}\ \psi(\ell_n)\ M(\C_s,\ell_n)}
           {\sum_n (n+1)^{-2}\ \psi(\ell_n)\ M(\C  ,\ell_n)}\  \; .
\ee 
with $\ell_n = 2^{-n}$.  Clearly, $t$ is a strictly increasing  
continuous function of $s$, and hence
defines a reparametrization of $C$. The denominator  satisfies
\be
\sum_n (n+1)^{-2} \ \psi(\ell_n)\, M(\C,\ell_n)
\ \le\ \sum_n (n+1)^{-2} \ < \ 2   
\ee
by the assumption~(\ref{ineq:M(C,l)}).
Consider two points  $\gamma(s_1)$ and $\gamma(s_2)$ 
(with $s_1<s_2)$ that are
at least $\Delta q$ apart, and let $\Delta t$ be the corresponding 
time difference.  For large  $n$ we observe that
\be
\ell_n\ <\ \Delta q \qquad   \Longrightarrow   \qquad 
 M(\C_{s_2},\ell_n)-M(\C_{s_1},\ell_n)\ge 1  \;  .
\ee
It follows that 
\bea
\Delta t\ & \ge\ & {1 \over 2} \sum_{n:\  \ell_n < \Delta q} 
(n+1)^{-2}\ \psi(\ell_n) \nonumber \\
 & \ge & \ { \psi(\Delta q/2) \over 
 2(\log_2 (4/\Delta q))^2 }\; ,
\eea
as claimed in \eq{ineq:Delta-gamma-2}. 
\hfill\blackbox

\masubsect{Tortuosity and box dimension}

In view of Theorem~\ref{thm:alphatau} it is important for us to 
have conditions implying finiteness of the tortuosity exponent.  
It is also of interest to have efficient estimates of the exponent's 
value.  Both goals are accomplished here through a criterion 
for the equality of $\tau(\C)$
with the upper box dimension $\dim({\C})$, which is relatively easier 
to estimate (and never exceeds $d$).  Some criterion is needed since 
in general the tortuosity exponent may exceed the 
upper box dimension, and may even be infinite~\cite{DDC}.  

\bigskip\noindent{\bf Definition} {\it \begin{itemize}
\item[i.]  We say that a curve $\C$  in $R^d$
exhibits a {\em $k$-fold crossing of power $\epsilon$, at the
scale $r \le 1$} if it traverses $k$ times some spherical shell
$D(x; r^{1+\eps}, r)$ (in the notation of \eq{eq:def-shell}).

\item[ii.]  A curve has the {\em tempered-crossing property}
if for every $0<\epsilon <1$  there are $k(\epsilon) < \infty$
and $0 < r_{o}(\epsilon) < 1$ 
such that on  scales smaller than $r_{o}(\epsilon)$
it has no $k(\eps)$-fold crossing of power $\epsilon$.
\end{itemize}
}

Note that the condition places restrictions on crossings at  
arbitrarily small scales; however, it is less restrictive
at smaller scales since it rules out only 
crossings of  spherical shells with increasingly large aspect ratio.

\begin{thm}  If a curve $\C$ has the tempered-crossing property, then
\be
\tau(\C) \ =\  \dim\C\  \; .
\label{eq:taubox}
\ee
In particular,  $\C$ admits H\" older continuous parametrizations
with every exponent $\alpha<(\dim(\C))^{-1}$.
\label{thm:taubox}
\end{thm}

\noindent {\bf Proof: }  Since $M(\C,\ell) \le N(\C,\ell)$, it is always 
true that $\tau(\C) \ \ge\  \dim\C\ $ (equation (\ref 
{ineq:dim-trivial})).    
To establish the opposite inequality we first prove that 
if a curve $\C$ has no $k$-fold crossings of power $\eps$ 
at the scale $\ell$ then 
\be
M(\C,2\ell)\ \le\ k\, N(\C,\ell^{1+\eps}) \; .
\label{eq:taubox-MN}
\ee  

To prove \eq{eq:taubox-MN} we recursively  partition the curve into
segments of diameter at most $2\ell$. The segments are defined
by a sequence of points $x_i$ along the curve.
We start with $x_1 = \gamma(0)$.  
After $x_1,\dots ,x_n$ are determined, the next point $x_{n+1}$  
is taken as the site of the first exit, after $n_n$,  of $\gamma$ 
from the  ball of radius $\ell$ about $x_n$; 
if $\gamma$ does not leave this ball, we terminate.

The number of stopping points produced by this algorithm  is clearly
an upper bound for $M(\C,2\ell)$.  In order to estimate 
this number, let us 
consider a covering of ${\C}$ by balls of diameter  
$\ell^{1+\eps}$.  Since there are no $k$-fold crossings of power
$\eps$ at the scale $\ell$, no such ball will contain more than $k$
of the stopping sites.  Hence \eq{eq:taubox-MN}. 

By the definition of the upper box dimension, for each 
$s>\dim\C$, the number $N(\C,\ell)$ of elements 
in a minimal covering  satisfies
\be
N(\C,\ell) \ \le\ K_s\ \ell^{-s} \  .
\ee
for some constant $K_s$ which depends on the curve. 
Therefore, for any $s >  \dim(\C)$
  \be
  M({\C},2\ell) \ 
 \ \le\ k\,K_s\,\ell^{-s(1+\eps)} \; ,
  \label{eq:taubox-singlescale}
 \ee 
 with some $K_s(\C)<\infty$.
Our assumptions imply that the exponent $s(1+\eps)$ 
 can be made arbitrarily close to 
$\dim \C$, and therefore $\tau(\C) \le \dim \C$.  That 
concludes the proof of \eq{eq:taubox}.  
The assertion about the H\"older regularity follows  from
Theorem~\ref{thm:alphatau}.
\hfill \blackbox

\bigskip\noindent{\bf Remark:\/} The proof of Theorem~\ref{thm:taubox}
shows that the tortuosity exponent can be bounded by the box dimension
under the weaker assumption  that for some integer $k$ and 
$\eps>0$, the curve
has no $k$-fold crossings of  power $\eps$ below some scale $r_o$.
In this case inequality~ (\ref{eq:taubox-singlescale}) implies the 
bound 
\be
\dim\C\ \le\ \tau(\C) \ \le \  (1+\eps)\  \dim\C   \; .
\label{eq:taubox-eps}
\ee

\masect{Regularity for curves in random systems}
\label{sect:random}

We now extend the discussion from a single curve to 
systems of random curves, in the terminology presented in the 
introduction.  Our first goal is to prove 
Theorem~\ref{thm:regularity}.  Following that 
we discuss the concept of the backbone, and thus improve 
the bounds on the dimension of curves. 

\masubsect{Proof of the main regularity result} 

An essential step   towards establishing regularity of random curves
consists of showing that under the hypothesis {\bf H1}, $k$-fold
crossings of spherical shells are rare in a sense which 
provides a probabilistic version of the tempered crossing condition. 
For this purpose, let us define the random variables
\be  
r_{\eps,k}^{(\delta)}(\omega)\ := \ 
\left\{ 
  \begin{array}{l}
	\inf \left \{ 0 < r \le 1  \ \left| \ 
	\begin{array}{c}
	   \mbox{ some shell $D(x; r^{1+\eps}, r)$, $x\in \Lambda$,
	    is traversed} \\
	   \mbox{ by $k$ distinct segments of curves in ${\cal  
	                               F_{\delta}(\omega)}$  } 
	\end{array} \right. \right\}  \\  
	\mbox{   } \\ 
	1 \qquad \qquad \mbox{\em if no such $k$ crossing occurs} 
  \end{array} 
 \right.  
\ee 
\begin{lem}  Let $\F$ be a system of random curves 
with variable short distance cutoff,
in a compact region $\Lambda\subset\R^d$ . 
Let $\eps > 0$, and assume that the condition (\ref {eq:suff-upper-1}) 
of \/ {\bf H1} holds for some $k<\infty$, and $\lambda(k)$ large enough 
so that  
\be \eps\,\lambda(k)-d > 0 \; .  
\label{eq:eps-k} 
\ee 
Then the random variable 
$r_{\eps,k}(\omega)$ is stochastically bounded away from zero, 
with 
\be
Prob_\delta \left( r_{\eps,k}^{(\delta)}(\omega) \ 
\le u\right)\ \le\ \ Const.(\eps,k)
\ u^{\eps\,\lambda(k)-d}\quad .
\label{eq:crossings}
\ee
\label{lem:crossings} 
\end{lem}

\noindent{\bf Proof:}\  We need to estimate the probability 
that there is a $k$-fold crossing of power $\eps$ at some
scale $ r \le u$.  Any such crossing gives rise to 
a crossing in a smaller spherical shell with discretized coordinates:
$D(x;3r_n^{1+\eps},r_n/2)$ with: $r_{n}= 2^{-n}$, 
$x \in (2r_n^{1+\eps}/\sqrt{d})\ \Z^{d}$ (where $\Z^d$ is the integer 
lattice in $\R^d$),   
and $n$ chosen so that $r_n < r \le r_{n+1}$. 
Using the assumption {\bf H1} and adding the probabilities over 
the possible placements of the discretized shells, we find:  
\bea
\P_{\delta}\left( \begin{array}{c}
\mbox{$\F_\delta(\omega)$  exhibits a $(k,\eps)$ 
crossing} \\ 
\mbox{at some scale $r \in (r_n,r_{n+1}]$}
\end{array} \right) 
&\le& \left(\frac{\sqrt{d}}{2r_n^{1+\eps}}\right)^{d} \ K_k \left( 
\frac{3r_n^{1+\eps}}{r_n/2}\right)^{\lambda(k)} \nonumber \\
&\le& Const. \ r_n^{\eps\,\lambda(k)-(1+\eps)d} \;  ,
\eea
where the constant depends only on $k$, $\lambda(k)$, and the 
dimension.  This bound bound decays exponentially in $n$. 
Its sum over scales  $r_n$ ($\delta\le r_n \le u$)
yields the claim.
\hfill\blackbox 

\bigskip\noindent 
{\bf Proof of Theorem~\ref{thm:regularity}}  
First let us note that the statement to be proven 
can be reformulated as follows. 

\noindent {\em Let $\F$ be a system of random curves in a compact region 
$\Lambda \subset \R^d$, with variable short-distance cutoff 
$\delta>0$, and assume the hypothesis $\bf H1$ is satisfied.  
Then for any $\eps > 0$ there is a random variable 
$\tilde \kappa_{\eps;\delta}(\omega)$, 
which stays stochastically bounded as $\delta \to 0$, with which 
the following tortuosity bound
applies simultaneously to all the curves $\C \in \F_{\delta}(\omega)$:
\be
M(\C,\ell) \ 
 \le \tilde \kappa_{\eps;\delta}(\omega) \ 
 {1 \over (\mbox{diam\/}\C)^{\lambda(1)+\eps}   }
  \ \ell^{-[d - \lambda(1) +\eps]}\; .
\label{eq:tortuosity-gen}
\ee    }

In this formulation, the H\"older continuity 
estimate of \eq{eq:holder-gen} is replaced by a tortuosity bound.  
The equivalence is based on Lemma~\ref{lem:uniform}, with  the function
$\psi(\ell)=Const.\  \ell^{\alpha}$ for which the inverse 
function is the power law with $s=1/\alpha$.  The 
logarithmic correction in \eq{eq:uniform-log} is absorbed 
through the ``infinitesimal slack'' we have in in the power law.    

Let now $\eps > 0$.  By the hypothesis {\bf H1}, 
there exists $k$ large enough so \eq{eq:eps-k} 
is satisfied.   
For such value of $k$, we learn from Lemma~\ref{lem:crossings} 
and Theorem~\ref{thm:taubox} (more specifically,  
\eq{eq:taubox-MN} there) that 
for $\ell$ small enough, i.e., $\ell < r_{\eps,k}(\omega)$, 
\be 
M(\C,2 \ell) \ \le \ k N(\C, \ell^{1+\eps})  \;  .
\label{eq:M1} 
\ee 
In the complementary range, $\ell\ge r_{\eps,k}(\omega)$, we use that
\be
M(\C,2\ell) \ \le\  M\left(\C,2 r_{\eps,k}(\omega)\right)\ 
\ \le\  k N(\C,r_{\eps,k}(\omega)^{1+\eps})\ .
\ee
It follows that
\be 
M(\C,2 \ell) \ \le \ A_{\eps,k;\delta}(\omega)\  
N(\C, \ell^{1+\eps}) \; ,
\label{eq:M2} 
\ee 
where the random variable 
\be 
A_{\eps,k;\delta}(\omega) \ = \ 
\left( \frac{\ell}{r_{\eps,k}(\omega)}\right)^{(1+\eps)d}
\ee
remains stochastically bounded as 
$\delta \to 0$ by Lemma~\ref{lem:crossings}.  

We shall now introduce some useful random variables 
which will permit to extract from \eq{eq:M2}  bound s
valid simultaneously for all curves 
$\C\in \F_{\delta}(\omega)$.   
Referring to the standard grid partition 
of $\Lambda$, let
\be 
\tilde N_{\delta}(r,\ell;\omega)\ := \ \left\{ 
 \begin{array}{l}
 	\mbox{the number of cubes $B$ of diameter $\ell$}   \\
 	\mbox{which meet a curve $\C \in  \F_{\delta}(\omega)$ with}\\ 
 	\mbox{diameter $diam(\C) \ge r$}
 \end{array}   \right.  
\ee 
Its expectation value is $\E( N_{\delta}(r,\ell) )$.  
Summing over scales $r_n\ge\ell_m\ge\delta$, with $r_n = \ell_n=2^{-n}$,
we define 
\be
U_{\delta}(\omega) \ := \ \sum_{m\le n}
{ \tilde N_{\delta}(r_n,\ell_m;\omega) \over \E( N_{\delta}(r_n,\ell_m) ) } 
\ (n+1)^{-2} (m+1)^{-2} \; .
\label{eq:defU}
\ee 
This random variable  stays stochastically bounded as $\delta \to 0$, 
by the Chebysheff inequality and the observation that the mean is 
independent of $\delta$:
\be 
\E( U_{\delta}) \ \le \ \left[ \sum_{n=1}^{\infty} 1/n^2 \right]^2 \; .
\ee 

For  the mean value of $\tilde N_{\delta}(r,\ell; \omega) $ 
we find
\bea
\E_\delta \left( \tilde N_{\delta}(r,\ell;\omega) \right) 
&= & \sum_{B\subset \Lambda; \ diam(B)=\ell}
 \P_\delta \left( 
 \begin{array}{c}  
 	\mbox{$B$ meets a curve $\C \in \F$}  \\
 	\mbox{with $diam(\C) \ge r$ } 
  \end{array}   \right)     \nonumber   \\
&\le&  \sum_{B\subset \Lambda; \ diam(B)=\ell} 
          K\,\left(\frac{\ell/2}{r/2}\right)^{\lambda(1)}\nonumber \\  
&\le& {K' |\Lambda| \over r^{\lambda(1)}}  
\left({1 \over \ell }\right)^{d-\lambda(1)} \;  .
\label{eq:meanN}
\eea

We now return to \eq{eq:M2}.   For curves with $diam(\F) \ge r$ we 
use
\be
N(\C, \ell) \ \le \ \tilde N_{\delta}(r,\ell;\omega) \ \le \ 
 \left(\log_2 {2/\ell)}\right)^2 \left(\log_2 {(2/r)}\right)^2  
     U_{\delta}(\omega) \ \E(  \tilde N_{\delta}(r,\ell) ) \; .
\label{eq:N-bound}
\ee 
(the last inequality based on the definition (\ref{eq:defU}). 
Combining the equations (\ref{eq:M2}), (\ref{eq:meanN}), and  
(\ref{eq:N-bound}), we learn:
\bea 
M(\C,2 \ell) \ & \le & \ 
\left[  (1+\eps)^2  
   \ K' \ |\Lambda| \
A_{\eps;\delta}(\omega)\ U_{\delta}(\omega) \right]  \ \times 
\nonumber \\ 
& & \qquad \times\ { \left( \log_2 {(2/\ell)}\right)^2  
       \left(\log_2 {(2/r})\right)^2 \over r^{\lambda(1)} } \ 
\left( {1 \over \ell } \right)^{(1+\eps)[d - \lambda(1)]} 
 \;  . 
\label{eq:almostthere} 
\eea 
The product of stochastically bounded variables is stochastically 
bounded, and the logarithm can be absorbed by adjusting $\eps$.  
Hence \eq{eq:almostthere} implies the claimed 
\eq{eq:tortuosity-gen}. \phantom{r}       \hfill\blackbox

\masubsect{Tortuosity of random systems and the backbone dimension} 

To summarize some of the results in a compact form, it may be 
useful to extend the notions of tortuosity and 
dimensions to systems of random curves with varying cutoff.  

\noindent {\bf Definition} {\em  For a system $\F$ of curves   
in a compact set $\Lambda \subset \R^d$ 
\begin{itemize}
\item[i.] The {\em upper tortuosity exponent  $\overline{\tau}(\F)$} 
is the infimum of $s>0$ for which the  random  variables 
\be 
 \sup \left\{  M(\C,\ell) \, \ell^s  \ | 
 	 \C \in \F_{\delta}(\omega)\; , \mbox{diam\/}(\C) \ge r \ \right\} 
\label{eq:def-taubar}
\ee 
remain stochastically bounded, as $\delta \to 0$ at fixed
 $ 0 < r \le 1$.  
\item[ii.] Similarly, the {\em upper box dimension} $\dim(\F)$ 
is defined through the boundedness of the 
variables given by
\be 
 \sup \left\{  N(\C,\ell) \, \ell^s  \ | 
 	 \C \in \F_{\delta}(\omega)\; ,  \mbox{diam\/}\C \ge r \  \right\}  
\ee 
with  $s,r$ as above.
\end{itemize}   }

The analysis carried above implies that if hypothesis {\bf H1} holds  
then the upper tortuosity exponent $\overline \tau(\F)$ is finite, 
and furthermore
\be
\overline {\tau}(\F) \ = \ \dim(\F)\ \le \ d - \lambda(1).
\label{eq:lambda1}
\ee 

The dimension estimate \eq{eq:lambda1}
reflects the fact that each point
on a curve $\C \in \F$  is connected a macroscopic distance away.  
It might seem that most  points on a curve are in fact at the end 
points of {\em two} line segments of macroscopic length.  This suggests 
an improved upper bound, in which $d-\lambda(1)$ is replaced by 
the smaller $d-\lambda(2)$.  However, one has to proceed
here with caution.  

There are two reasons for which $d-\lambda(2)$ may 
not provide a valid upper bound for the dimension:  
\begin{itemize}  
\item[i.] 
The union  $\bigcup_{\C\in \F_\delta(\omega);\ \mbox{diam\/}\C > r} \C $
 may be dominated by  
the collection of the {\em endpoints} of curves in $\F$ which are only 
{\em singly connected} a macroscopic distance away.  
For instance, that would occur if the connected clusters to which the 
curves of $\F_\delta(\omega)$ are restricted have many short branches 
(one could call this the {\em broccoli effect}).  
\item[ii.]  Certain curves $\C \in \F$ may be rougher 
at the their ends, where only one segment is accommodated in the 
available space, than in their interior. We expect this to be the case for 
some examples of self-avoiding paths.  
When that happens it will not be true that ``most of the curve'', 
as counted by covering boxes, consists of its interior.
\end{itemize}

Nevertheless, the proposed bound is obviously valid for 
the union of the interior parts of the curves, if that is defined 
as the collection of point on $\C$ whose distance to the endpoints is 
at least some $a>0$, which remains fixed  as $\delta \to 0$.  
The proof is by a direct adaptation of argument used in the 
proof of Theorem~\ref{thm:regularity}, making the 
suitable correction in \eq{eq:meanN}.

A situation like that has been addressed in the percolation context 
through the concept of the  {\em backbone}.  The term is used  to
distinguish between a {\em spanning cluster}, i.e., a cluster 
connecting two opposite faces of a macroscopic size cube which 
typically contains many {\em dangling ends}, and the smaller set of 
bonds which carry current between the faces~(\cite{Sta77}).

A mathematically appealing formulation is possible in the continuum 
limit  (at $\delta=0$), for which we define the backbone 
${\cal B}(\omega)$ of the system of curves  $\F(\omega)$ as the 
union of all interior segments of curves $\C\in \F(\omega)$.   

For the backbone, the Hausdorff and box dimensions need not 
coincide.  Since the statement is closely related with the 
considerations of this section we present it here, even though it 
anticipates the construction which is better described in the next 
section.
\begin{thm}
In the scaling limit [defined in the next section] 
\be
dim_{\cal H}{\cal B}(\omega) \ \le \ d - \lambda(2) \quad (a.s.) 
\; ,
\label{eq:bbhausdorff} 
\ee
whereas 
\be 
 \dim{{\cal B}(\omega) } \ = \ \dim{{\cal F(\omega) } } \ \le  \ 
 d - \lambda(1) \quad (a.s.) .
\label{eq:bbbox}
\ee 
\label{thm:bb}  
\end{thm} 
The last inequality can be saturated.   

\noindent {\bf Proof:} 
Equation (\ref{eq:bbhausdorff}) follows from the continuity of the
Hausdorff dimension under countable unions, and the 
the previous observation on  the dimension of the sets defined 
with fixed macroscopic cutoffs.  Equation (\ref{eq:bbbox}) holds 
since the box dimension of a set equals that of 
its closure, which for ${\cal B(F(\omega))}$ is the union of
all curves in ${\cal F}$.   
\hfill \blackbox

\masect{Compactness, tightness, and scaling limits}
\label{sect:compact}

We now turn to the construction of scaling limits for a 
random system of curves.  Such a system  is described by a collection of 
probability measures $\mu_\delta$ on the space 
of configurations of curves, $\Omega_\Lambda$ defined in the 
introduction.   We shall see that the tortuosity 
bound~\ref{eq:tortuosity-gen}
derived in Theorem~\ref{thm:regularity} 
allows one to conclude the existence of limits for $\mu_\delta$.  

The key to the proof of Theorem~\ref{thm:scaling} is the relation 
of  the space of curves with the space of continuous functions,  
$C([0,1],\Lambda )$, 
and the well developed theory of probability measures on the space 
of closed subsets of a complete separable metric space.  We recall 
some of this theory below.  The first step is 
the following counterpart to Arzel\`a-Ascoli theorem. 
 
\begin{lem} {\rm (Compactness in ${\cal S}_\Lambda$)} 
A closed subset ${\cal K} \subset {\cal S}_\Lambda$, of the 
space of  curves in a compact $\Lambda\subset\R^d$, is compact 
if and only if there exists a function  
$\psi:(0,1]\to(0,1]$ so that for all $\C \in {\cal K}$,
\be 
M(\C,\ell)\ \le\ \frac{1}{\psi(\ell)}\quad 
\mbox{for all}\ 0<\ell\le 1\; .
\label{eq:compact}
\ee
\label{lem:compact}
\end{lem}

\noindent{\bf Proof:} \  
We first show that if a closed set  ${\cal K}\subset 
{\cal S}_\Lambda$ consists of curves  
satisfying uniform tortuosity bounds then  ${\cal K}$ is compact.
It suffices here to show that each sequence of curves in ${\cal K}$ 
has an accumulation point in ${\cal S}_\Lambda$.  The limit will 
be in  ${\cal K}$ because  ${\cal K}$ is closed.

By Lemma~\ref{lem:uniform}, we can parametrize each curve in the 
sequence by a continuous function  satisfying 
the corresponding continuity condition \eq{ineq:Delta-gamma-2}. 
That yields an equicontinuous family of functions in 
$C([0,1],\Lambda)$.  Applying the Arzel\`a-Ascoli theorem we 
deduce the existence of a uniformly convergent subsequence.  
It is easy to see that the curves defined by these functions also  
converge, with respect to the metric on ${\cal S}$.   

In the converse direction (which we do not use in this work), 
we need to show that if ${\cal K}$ is compact 
then $M(\C,\ell)$ is uniformly bounded on it.  That follows from 
Lemma~\ref{lem:semicontinuous}, which shows that: 
i)  $M(\C,\ell)$ is bounded above by $\tilde M(\C,\ell/3)$, ii)
since $\tilde M(\C,\ell)$  is upper semicontinuous by 
Lemma~\ref{lem:semicontinuous}, it
achieves its  supremum on the compact set ${\cal K}$.
\hfill\blackbox

\bigskip  Standard arguments,  such as used 
for $C([0,1], \R^d)$, show that the space of curves 
${\cal S}_{\R^d}$ is a complete and separable metric space.  
The completeness and separability of ${\cal S}_{\R^d}$ are 
passed on to $\Omega_\Lambda$.   For this space we get the 
following characterization of compactness.

\begin{lem} {\rm (Compactness in $\Omega_\Lambda$)} 
A closed subset $\tilde A$ of $\Omega_{\Lambda}$  is compact, if and 
only if there exists some $\psi: (0,1] \to (0,1]$ for which 
each  configuration $\F\in\tilde A$
consists exclusively of curves satisfying  abound of
the form \eq{eq:compact}.  
\label{lem:compact-omega}
\end{lem}

\noindent{\bf Proof:} \  
The claim follows from  the basic property of the 
Hausdorff metric, under which the closed subsets
of a compact metric space form a  compact space.
\hfill\blackbox

The scaling limit 
we are interested in is taken in the space of probability measures 
on  $\Omega_{\Lambda}$, for compact ${\Lambda}\subset \R^d$.  
Our discussion will now make use of a number of useful general
concepts and results.  Let us just briefly list those.  
A thorough treatment can be found in ref.~\cite{Bill}.  

A family of probability measures $\{\mu_n\}$ is said to be 
{\em tight}  
if for every ${\epsilon}$  there exists a compact set $A$ so that
$\mu_n(A)\ge 1-\eps$.  

The sequence $\mu_n$ is said to {\em converge to} $\mu$ if 
$\lim_{n\to \infty} \int f \ d\mu_n = \int f \ d\mu$
for every continuous function $f : \Omega \to \R$.  
If the distance function is uniformly bounded, as is the case
for measures on $\Omega_{\Lambda}$ with compact $\Lambda$, 
this convergence statement is equivalent to the existence of a
coupling as described in the introduction, below 
the statement of Theorem~\ref{thm:scaling}.  

A collection of measures is said to be {\em relatively compact}
if every sequence has a convergent 
subsequence.  Tightness and compactness are equivalent in this 
general setting:

\noindent {\bf Theorem} (Prohorov \cite{Proh}, see~\cite{Bill})\ 
{\em A family of
probability measures on a complete separable metric 
space is relatively compact if and only if it is tight.}
\bigskip

Thus, in order to prove Theorem~\ref{thm:scaling}, we need to show 
that for each $\eps > 0$, up to remainders of probability $\le \eps$ 
the measures $\{\mu_\delta\}$ are  supported on a common compact subset 
(of the space of configuration),  which may depend on $\eps$.

\noindent{\bf Proof of Theorem~\ref{thm:scaling} } By  
Theorem~\ref{thm:regularity} and point {\em ii.} of
the remark
following it, there is for each $s>d$ and $\eps>0$ a choice of $K<\infty$, 
such that all curves in the random configuration
$\F(\omega)$ drawn with the probability measure 
$\mu_{\delta}$  can  be parametrized  H\"older continuously with
exponent $s$ and H\"older constant $\kappa_{s;\delta}$ ,
as in \eq{eq:holder-gen}.
By Lemma~\ref{lem:uniform}, this implies that
\be 
M(\C,\ell)\ \le\ K_r \ell^{-s}\ \mbox{\rm for all curves}\ 
\C\in\F_\delta(\omega) \; ,
\label{eq:Ks} 
\ee 
except for a collection of configurations 
whose total probability is $\le \eps$.  By Lemma~\ref{lem:compact-omega},
the set $A_\Lambda(K,s) \subset \Omega_{\Lambda}$ of  
all configurations consisting only of curves that satisfy  \eq{eq:Ks}. 
is compact. In other words,
finite upper tortuosity of $\F$ implies that
upon truncation of small remainders the measures 
$\mu_{\delta}$ are  supported in the compact 
sets of the form $ A_{\Lambda}(K,s) $.  
(Note that $K<\infty$ needs to be 
adjusted depending on $\eps$ and the choice of $s$.)  
This proves that the family $\mu_{\delta}$ is tight. By Prohorov's
theorem, that is equivalent to compactness. 

To see that the limiting measure is supported on curves that can be 
parametrized  H\" older continuously with any exponent less
than $1/(d-\lambda(1))$, consider the collections $\F^{(r)}$ of
curves of diameter a least $r$. The above argument shows that
the measure restricted to this collection is almost supported
on $A_\Lambda(K(r),s)$ for any $s>d-\lambda(1)$, and $K(r)$ large enough.
By Prohorov's theorem, the limiting measure is supported
on $\bigcup_{K>0} A_\Lambda(K,s)$, which proves the claim by
Lemma~\ref{lem:uniform}.
\hfill\blackbox

\bigskip Let us remark that the notion of convergence  we use here
(technically it is called {\em weak convergence} on the space of
measures on $\Omega_\Lambda$) is  quite strong, due to our choice of
topology on $\Omega_\Lambda$.  
As \eq{eq:vasershtein} makes it clear, for $n$ large typical 
configurations of $\F_{\delta_n}$ are close to typical configurations
of the scaling limit -- close in the sense of the Hausdorff
metric induced on the space of configurations, $\Omega_\Lambda$,
by the uniform metric in the space of curves $S_\Lambda$. 
This sense of convergence is stronger than that defined through 
the joint probability distributions of  finite collections 
of macroscopic crossing events.  In this respect,  
the notion of convergence used here 
is  reminiscent of the sense in which Brownian motion is 
proven to  approximate random walks, in Donsker's theorem~\cite{Donsker}.

\masect{Lower bounds for the Hausdorff dimension of curves}
\label{sect:sparse}

Our next goal (the third theme of this work) is to prove 
the statement of Theorem~\ref{thm:roughness},
that in a system satisfying 
the hypothesis {\bf H2}, almost surely none of the curves which 
appear in the scaling limit are of Hausdorff dimension lower 
than some $d_{min}>1$.

The proof is split into two parts. The first part, 
carried out in this section, consists of  
measure-theoretic analysis based on the assumption that a 
certain auxiliary deterministic 
condition is satisfied for a given curve.  
In the next section the proof is completed 
with a probabilistic argument showing that in a system of
random curves satisfying
the hypothesis {\bf H2}, the auxiliary 
condition is met almost surely.  

\masubsect{Straight runs}

Standard examples of curves of dimension greater than one are 
curves whose segments deviate from straight lines
proportionally on all scales.
However for random systems (and other setups) that 
criterion is too restrictive since one may expect  exceptions 
to any rule to occur on many scales.    
The criterion which we develop here is the 
{\em sparsity of straight runs}, which is an abbreviated expression 
for the absence of sequences of nested straight runs occurring over 
an excessively dense collection of scales. The concept is defined with a 
macroscopic scale $L>0$ and a shrinkage 
factor $\gamma>1$, used to  specify a sequence of length scales:
\be
L_k \ = \  \gamma^{-k}\,L_o  \; ,
\ee  
and an integer $k_o$, used to allow exceptions above a certain
scale.

\bigskip\noindent {\bf Definition} \ {\it A curve in $\R^d$  
is said to exhibit a  {\em straight run} 
at scale $L$ ($=L_k$ for some $k$), 
if it traverses some cylinder of length $L$ and cross sectional
diameter $(9/\sqrt{\gamma})\,L$, in the 
``length'' direction, joining the centers
of the corresponding faces. 
Two straight runs are {\em nested} if one of the defining  
cylinders contains the other.

We say that straight runs are 
$(\gamma,k_o)$-{\em sparse}, down to the 
scale $\ell $, if $\C$ does not exhibit any nested collection  
of straight runs on a sequence of scales 
$L_{k_{1}}> \ldots > L_{k_n}$, with   $L_{k_n}\ge \ell$ and 
\be  
n \ge \frac{1}{2}  \max\{ k_n, k_o\} \; \ .
\ee 
}
\bigskip

Following is the deterministic result, which 
is stated here only  in the continuum ($\delta=0$).   For systems
of random curves we will make use of the more detailed information
which appears in the proof (see \eq{eq:lower-cap}).

\begin{thm}
If the straight runs of a given curve $\C$ are
$(\gamma,k_o)$-sparse , then
$dim_{\cal H}\C \ge s $, with $s$ given by
\be
 \gamma^s =
\sqrt{m(m+1)}\; ,
\ee
and $m$ an integer strictly smaller than $\gamma$.
\label{thm:lower-det}
 \end{thm}
 
Clearly, if for some integer $m$ the above condition is 
met for all $\gamma > m$  then the bound becomes 
\be
dim_{\cal H}\C \ \ge 
 \ 1 + \frac{ \ln (1+1/m)}{2\, \ln m}  \;  \ .
\ee 

We will prove Theorem~\ref{thm:lower-det}
by cutting the given curve $\C$ into a hierarchical family of
subsegments at different scales, with segments at the same
scale separated by a certain minimal distance.
This family defines a Cantor-like (i.e., closed, perfect, 
and totally disconnected) subset $\tilde C$ of $\C$.
If $\C$ contains no straight runs at all,
a scaling argument easily shows that the dimension  of $\tilde C$, and
hence the dimension of $\C$, exceeds one. We use capacity arguments
to show that this holds also under the weaker condition that
straight runs are sparse. For the construction of
the family of  subsegments which defines the fractal subset
$\tilde \C$, we modify the exit-point algorithm from the proof 
of Theorem~\ref{thm:taubox}.

\masubsect{Construction of fractal subsets}  

Let  $\gamma$ be a positive number,  $m$  an integer in 
$[\gamma/2, \gamma ]$, and $k_{max}$ a positive 
integer.  By an iterative procedure we shall construct for a given curve 
$\C$  a nested sequence $\Gamma_o,\ldots ,\Gamma_{k_{max}}$ of 
collections of segments of $\C$, at scales
\be 
 L_k \ = \ \gamma^{-k}\,L_o\; ,\quad 
  k = 0,\dots, k_{max} 
 \ ,
\ee 
with $L_o = {\rm diam}(\C)$, having the  following properties.
\begin{itemize}
\item[i.]
Each $\Gamma_k$ is a collection of segments of diameter
at least $L_k$.
\item[ii.]
In each generation (as defined by $k$),  distinct segments are at 
distances at least $\eps L_k$ with $\eps = \gamma/m - 1$.
\item[iii.] 
Each segment of $\Gamma_k$ ($k>1$) is contained in one of the 
segments of $\Gamma_{k-1}$.  The number of immediate 
{\em descendants} thus contained in a given element  of $\Gamma_{k-1}$ 
is at least $m$, and very frequently at least $(m+1)$.
\end{itemize}

To define $ \tilde{\C}$ let $\C_k$ be the union of the segments 
in $\Gamma_k$.  Then $ \tilde{\C}=\cap_{k\le k_{max}} \Gamma_k $. 
In the construction, we find it convenient to use
the {\em span} of a curve, which we define
to be the distance between the curve's end points, in place of the diameter.

\begin{lem}  {\rm (Construction of  $\tilde{\C}$)} \ 
There is an algorithmic construction which for each curve 
$\C$, yields a sequence of 
collections of segments with the above properties {\em i. - iii.\thinspace}, 
and with the further property that unless a segment 
$\eta \in \Gamma_k$ exhibits a straight run of 
scale $L_k$ the number of its descendants is
at least $m+1$.
\label{lem:fractal}
\end{lem}

\noindent {\bf Proof:}  We may assume (by trimming) that the span of 
$\C$ equals $L_o$.  Let $\Gamma_o$  consist of only 
one element: a segment which starts at one end of the curve and 
stops upon the first exit from a ball of radius ${\rm diam}(\C)$.
Once $\Gamma_k$ has been constructed, we form 
$\Gamma_{k+1}$ by selecting for each 
element $\eta \in \Gamma_k$ a collection of  descendants
$\eta_1,\dots \eta_N$, which are subsegments of $\eta$ cut by 
two sequences of points $x_j$ and $y_j$, strung along it 
in the order:  $y_1 < x_1 < y_2 < x_2 < \ldots$.  
The cutting points are selected by the following procedure.

 \begin{figure}[t]
      \begin{center} 
      \leavevmode
  	\epsfysize=2in
      \epsfbox{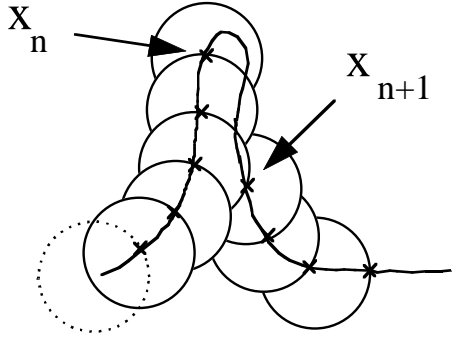}
  \caption{\footnotesize  The algorithm used 
  for marking the points $x_j$ in the construction of $\Gamma_k$. }
  \label{fig:algorithm}
  \end{center}
  \end{figure}  
    
We let $y_1$ be the starting point, and $x_1$ the first exit 
of $\eta$ from the ball of radius $L_{k+1}=L_k/\gamma$ centered  
at $y_1$.  Then, recursively, we chose 
$x_n$ as the first point on $\eta$ having distance
at least $L_k/m$ from the already constructed subsegments 
$\eta_1,\dots ,\eta_{n-1}$; terminating if  no such point can be  
found. The point $y_n$ is selected as the last entrance, 
prior to $x_n$ into the ball of radius $L_{k+1}$ 
centered at $x_n$.

 \begin{figure}[t]
      \begin{center} 
      \leavevmode
  	\epsfysize=2in
      \epsfbox{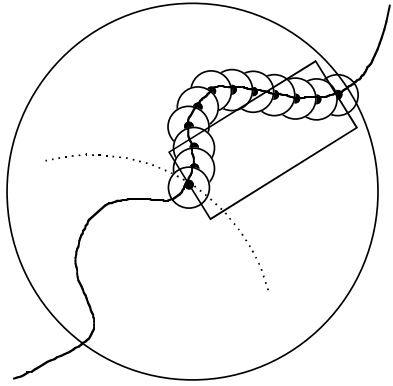}
  \caption{\footnotesize  The subdivision of an element of 
  $\Gamma_k$ into $\Gamma_{k+1}$.  Unless there is
  a straight run in a cylinder positioned as indicated, the number
  of elements increases at a higher rate than the factor by 
  which the radius shrinks ($\gamma$).  Under the hypothesis {\bf H2} 
  straight runs are  sparse in a sense which 
  permits to derive  a lower bound  on the Hausdorff dimension  based 
  on this picture. }
  \label{fig:fsubdivision}
  \end{center}
  \end{figure}  
  
It can be verified that the sequence of subsegments 
$\eta_j$, $j=1,\ldots,N$, with the endpoints $\{x_j, y_j \}$ 
have the properties:  
\begin{itemize}
\item[i.]
each $\eta_N$ spans  distance  $L_{k+1}$,
\item [ii.] 
the distance of each point on $\eta_n$ to the union of the 
segments $\eta_1\dots \eta_{n-1}$ lies between $\eps L_{k+1}$ 
and $L_k/m$, while
\item[iii.] the distance from $y_n$
to the starting point $x$ is at most $nL_{k-1}/m$.
\end{itemize}
  
We need to estimate the number of segments generated by the above 
procedure.  It is easy to construct from the collection of segments
a polygon with step size at most $L_{k-1}/m$, connecting the 
endpoints of $\eta$.  Choose as the vertex before the last 
the point $x_N$, and for  any given vertex in $\eta_n\  (n>1)$, 
select the preceding one from  some $\eta_i$ with $i<n$  so that 
the resulting leg
has length at most $L_k/m$; if $n=1$, terminate and
use $y_1$ as the initial point. Clearly, the polygon
has at least $m$ interior vertices, and hence the number $N$ of 
subsegments is at least $m$.

Assume now that $\eta$ does not exhibit a straight run at scale
$L_k$. We claim that the number $N$ of descendants is at least $m+1$.
By construction, at least one of the segments $\eta_n$ 
has a distance of less than $L_k/m$ from the lateral boundary
of the cylinder of width $9/\sqrt{\gamma}$
defining a straight run from$x$ to $y$. If that segment
contributes a vertex to the polygon,
then this vertex must lie outside the cylinder
of width $(9/\sqrt{\gamma}-4/m)\,L_k\ \ge\ (2/\sqrt{m})\,L_k$.
Then the polygon has length at least $\sqrt{1+4/m}\,L_k\ge (1+1/m)\,L_k$,
and hence contains at least $m+1$ interior vertices, coming from
distinct  subsegments.
On the other hand, if some subsegment does not contribute a vertex to the
polygon, we also have $N>m$.  
This completes the proof of the lemma.
\hfill\blackbox

One may think of the elements of $\cup_{k} \Gamma_k$ as vertices of 
a tree, with 
the root in $\Gamma_o$, and edges joining each segment to its 
immediate descendants.   For any two points $x, y \in \tilde{\C}$ 
which are not in the same element of $\Gamma_{k_{max}}$
\be 
|x-y| \ \ge \ \eps L_{k(x,y)}  \;  ,
\ee 
where $k(x,y)$ is the index of the first generation at which 
the two points are separated.  
Following are two general results which we shall use to 
estimate the dimension 
of $\tilde{\C}$. 
 
\masubsect{Energy estimates}

For a metric space $A$ and $\ell\ge 0$, let $\mbox{Cov}_{\ell}(A)$ 
denote the collection of coverings of $A$ by sets of diameter not 
smaller than $\ell$. By the definition of the Hausdorff dimension, a lower 
bound on $dim_{\cal H}{A}$ means that 
for some $s > 0 $ the quantity 
\be
\inf_{\{B_j\} \in \mbox{Cov}_{\ell}(A)} \sum_{j} 
\left( diam B_j \right)^s 
\ee 
{\em does not tend to $0$} as $\ell \to 0$ (in which case 
$dim_{\cal H}{A}\ge  s $).  It is difficult to use
this definition directly to find lower bounds on the Hausdorff dimension.  
We shall therefore make use of the relation of Hausdorff measures 
with capacities and deduce a lower bound on dimension
from an upper bound on the {\em energy} of a judiciously chosen 
probability measure ({\em charge distribution}) 
supported on the set $A$.   

\begin{lem} 
For $s>0$ and $\ell\ge 0$, let the capacity  ${\CAP}_{s;\ell}(A)$
of a subset of $\R^d$ be defined by: 
\be
\frac{1}{\CAP_{s;\ell}(A)} \ =
 \ \inf_{\mu \ge 0, \int_A d\mu = 1} \int \int_{A\times A}
 \frac{\mu(dx) \mu(dy)}{\max \{|x-y|,\ell\}^s} 
\label{eq:cap-def}
\ee 
Then, for every collection of sets $\{B_j \}$ covering  $A$, 
with $\min_{j} \mbox{diam}(B_j) \ge \ell$:
\be 
\sum_{j} \left( \mbox{diam} B_j \right)^s \ \ge \ {\CAP}_{s;\ell}(A) \ .
\label{eq:cap-Hausdorff}
\ee 
\label{lem:cap-Hausdorff}
\end{lem} 

\noindent (The case $\ell=0$ can be found in Falconer~\cite{Falconer}.
The statement is related
to the theorem of Erd\H{o}s and  Gillis \cite{Erdos-Gillis} that 
the $s$-dimensional 
Hausdorff measure of $A$ is infinite whenever $\CAP_{s;0}$ is positive.)

\noindent {\bf Proof:} 
By monotonicity, it clearly suffices 
to prove  \eq{eq:cap-Hausdorff} for any covering by {\em disjoint} sets.  
Let $\{ B_j \}$ be such a collection, and $\mu $ a probability 
measure supported on $A$.  Then
\bea
 \int \int_{A\times A} \frac{\mu(dx) \mu(dy)}{\max \{|x-y|,\ell\}^s} 
  \ & \ge & \ \sum_j  \int \int_{x,y \in B_j}
 \frac{\mu(dx) \mu(dy)}{\max \{\mbox{diam}(B_j),\ell\}^s}   
 \nonumber  \\ 
   \ & = & \  \sum_j  \frac{\mu(B_j)^2}{\mbox{diam}(B_j)^s} \ . 
\label{eq:mu1}
 \eea  
We also have 
\be
1 \  =  \ \left(\sum_j \mu(B_j) \right)^2
 \le  \ \left(\sum_j 
\frac{\mu(B_j)^2}{\mbox{diam}(B_j)^{s}} \right) \ 
\left( \sum_j \mbox{diam}(B_j)^{s}  \right)
\; 
\label{eq:mu2}
\ee 
(by the Schwarz inequality).  
Combining the last two relations we learn that  
\be 
 \sum_j 
 \mbox{diam}(B_j)^{s} \ \cdot\ 
  \left(\int \int_{A\times A} 
  \frac{\mu(dx) \mu(dy)}{\max\{|x-y|,\ell\}^s}\right) \ \ge\ 1
\ee 
{\em for any probability measure supported on $A$, 
and any covering of $A$} by sets with diameters $\ge \ell$. 
Minimizing over $\mu$  one obtains the relation claimed in 
\eq{eq:cap-Hausdorff}. 
\hfill\blackbox

\bigskip 
\begin{lem}
Let $A$ be a compact subset of $\R^d$.  Assume there is a sequence 
$\Gamma_o, \ldots, \Gamma_{k_{max}} $ of (nonempty)
collections of 
closed disjoint  subsets of $\R^d$, 
such that for each $k=0,\ldots, k_{max} (\le \infty)$: 
\begin{itemize}
\item[i.]  Each element of  is contained in some element  
of $\Gamma_{k-1}$, and each element of $\Gamma_{k}$ contains at least 
one such ``descendant''.   
\item[ii.]  Any two distinct sets in  $\Gamma_{k}$ are a distance at 
        least $\eps L_k$ apart, where $L_k=\gamma^{-k},L_o$ with
some $L_o>0$, $\gamma > 1$ and $0<\eps \le\gamma$.
\item[iii.]  For each element $\eta \subset \Gamma_k$:  
$ \eta  \CAP A \neq  \emptyset$.
\end{itemize}
For points $x\in \cup_{\eta\in \Gamma_k} \eta $, let $n_k(x)$ be the 
number of immediate descendants of 
the set containing $x$ within  $\Gamma_{k-1}$.   
Assume, furthermore: 
\begin{itemize}
\item[iv.] 
there is some $\beta > 1$, such that 
\be 
\prod_{j = 1}^{k} n_j(x) \ge \beta^k  \; , \ \ \ \ 
\mbox{for all $k= k_o,\ldots, k_{max}$ }
\label{eq:product} 
\ee 
with some common $k_o$, 
whenever $x \in \cup_{\eta \in \Gamma_k} \eta $. 
\end{itemize}
Then, for $s>0$ such that $\gamma^s < \beta$, and 
$\ell = \gamma^{-k_{max}}L_o$:
\be
\CAP_{s;\ell}(A) \ \ge \ 
 (\eps\,L_o)^{s} \left[ \gamma^{s\,k_o}  \ + \ 
    \frac{ \beta}{1-\gamma^{s}/\beta}\right]^{-1} \; .
\label{eq:cap-bound}
\ee 
\label{lem:cap-bound}
\end{lem}

\bigskip\noindent{\bf Remark:\/}
It should be appreciated that $\ell$ and $k_{max}$ 
do not appear on the right side in  \eq{eq:cap-bound}.  
If straight runs are sparse on
all scales
(that is, $k_{max}=\infty$), then the limit $\ell \to 0$ of
\eq{eq:cap-Hausdorff}  yields a bound on
the s-dimensional Hausdorff measure of $A$.

\noindent{\bf Proof:} 
For a bound on the capacity it suffices 
to produce a single probability measure supported on $A$ 
with a correspondingly small  ``energy integral''   (see \eq{eq:cap-def}).
We construct the  measure $\mu$ so   
that for each $\eta \in \Gamma_k$ the total measure of
$A\cap \eta$ is distributed evenly among its immediate descendants.  
This means that for each $k=0,\ldots,k_{max}$ and 
each $\eta \in \Gamma_k$
\be 
\mu(\eta) \ = \ \prod_{j = 1}^{k} n_j(\eta)^{-1} \; ,
\label{eq:def-mu}
\ee
where (fpr $j\le k$), the number
$n_j(\eta)$ is the constant value which $n_j(x)$ takes for 
$x\in \eta$.
To specify the measure uniquely, we designate as its support 
$\{ x_{min}(\eta) \ | \ \eta \in \Gamma_{k_{max}} \}$, where 
for each $\eta\in\Gamma_{k_{max}}$ , the point $x_{min}(\eta)$  
is the earliest point in $\eta$, with respect  to
the lexicographic order of $\R^d$.

For $x,y \in A$, if the two points are in separate elements 
of $\Gamma_{k_{max}}$ we let $k(x,y)$ denote the index 
of the level at which they separated.  
In estimating the energy integral we shall use the bound:
\be
|x-y| \ \ge  \ \eps \ L_{k(x,y)} \; . 
\ee 
for points which are separated in $\Gamma_{k_{max}}$.  
Otherwise, we use $\max \{|x-y|,\ell\} \ge \ell$.  Thus 
\bea
 {\mathcal E}(\mu) &\equiv&\   \int \int_{A\times A} \ 
 \frac{\mu(dx) \mu(dy)}{\max \{|x-y|,\ell\}^s} \nonumber \\
  & \le & \  \int \int_{k(x,y)\le k_{max}}  
(\eps L_{k(x,y)})^{-s} \ \mu(dx) \mu(dy) \ + \  \nonumber \\ 
  & & \qquad \qquad +  \ \sum_{\eta \in \Gamma_{k_{max}} } 
  L_{k_{max}}^{-s}\ \int \int_{\eta \times \eta } 
  \mu(dx) \mu(dy)    \; . 
\eea  
Splitting the first   integral  on the right according to the 
value of  $k(x,y)$ (separating out the case $k(x,y)\le k_o$)
and replacing $L_k$ by $\gamma^{-k} L_o$ throughout, we obtain
\be
  {\mathcal E}(\mu) 
  \   \le   \ (\eps \,L_o)^{-s} \ \gamma^{s\,k_o}\ + 
  \ \sum_{k=k_o+1}^{k_{max}} 
   (\eps \,L_o)^{-s} \ \gamma^{s\,k} \sum_{\eta \in \Gamma_{k-1}}
  \mu(\eta )^2 \ + \ L_o^{-s}\,\gamma^{s\,k_{max}}
  \sum_{\eta  \in \Gamma_{k_{max}}}
  \mu(\eta )^2 \; .
\label{eq:energy-est} 
\ee  
Since $\eps\le \gamma$, the last term on the right hand
side of (\ref{eq:energy-est}) can be replaced there by adding 
the term $k=k_{max+1}$ to the preceding sum. 
Finally, we use the assumption \eq{eq:product} together with 
the definition of the measure in \eq{eq:def-mu} and
\be
\sum_{\eta \in \Gamma_k} \mu(\eta) \ = \ 1\;  
\ee 
to see that
\be 
\sum_{\eta\in\Gamma^k} \mu(\eta )^2 \ = \ 
\prod_{j = 0}^{k} n_j(\eta)^{-1} \,\sum_{\eta\in\Gamma_k} \mu(\eta)\  \le \  
\beta^{-k}\ .
\ee 
This yields a geometric series bound for the sum over $k$ in 
\eq{eq:energy-est}, which results in the bound stated in 
\eq{eq:cap-bound}. 
\hfill\blackbox

\noindent{\bf Proof of Theorem~\ref{thm:lower-det}:}
Let $\cal C$ be a curve where straight runs are 
$(\gamma,k_o)$ sparse down to scale 
$\ell=\gamma^{-k_{max}}$.
The hierarchical construction of
Lemma~\ref{lem:fractal} results in a fractal subset 
$\tilde {\cal C}$ of $\C$.
Since straight runs are sparse by assumption, $\tilde C$ satisfies
the branching condition (\ref{eq:product})
of Lemma~\ref{lem:cap-bound},  with the the value of $\beta$ 
defined by the relation: 
\be  
\beta \ = \ \sqrt{m(m+1)}\ .
\ee 
Thus, Lemma~\ref{lem:cap-bound} 
implies that for any $s$ such that $\gamma^s < \beta$:
\be 
\CAP_{s;\ell}({\C})  \ \ge  \  
(\eps\,L_o)^s\ \left[ \gamma^{s\,k_o}  \ + \ 
    \frac{ \beta}{1-\gamma^{s}/\beta} \right]^{-1} \; ; 
\label{eq:lower-cap}
\ee 
this inequality holds for all $\ell \in (\gamma^{-k_{max}},1]$.  
By Lemma~\ref{lem:cap-Hausdorff}, the same lower bound holds
for  $\inf_{\{B_j\} \in \mbox{Cov}_{\ell}(\tilde{\cal C})} \sum_{j} 
\left( diam B_j \right)^s $.  Since we may choose $k_{max}$ as 
large as we please, the
$s$-Hausdorff measure of $\C$ is positive, and hence,
the Hausdorff dimension is at least $s$.
\hfill\blackbox

\masect{Lower bounds on curve dimensions in random systems}
\label{sect:lower}

We shall now combine the previous deterministic results with a 
probabilistic estimate, and prove Theorem~\ref{thm:roughness}.
The proof consists of showing that, with high 
probability,
straight runs are sparse, and then applying the results of the 
previous section.
 
\begin{lem}  (Sparsity of straight runs.) 
Assume a system of random curves in a compact set $\Lambda\in\R^d$
satisfies the hypothesis {\bf H2}. 
For $\gamma>4d$, define
a sequence of length scales $L_k=\gamma^k$. Then
there are constants $K_\Lambda, K_1< \infty, K_2>0$, with which for 
any fixed sequence $k_1< k_2<\ldots< k_n$ 
\be
Prob_{\delta}\left(  
\begin{array}{c}
	\mbox{there is a nested sequence} \\ 
\mbox{of straight runs at scales $L_{k_1},\ldots, L_{k_n} $ }
\end{array}  \right) \ \le \ K_\Lambda\ \gamma^{2d\, k_n}\, 
e^{(K_1-K_2\sqrt{\gamma})\, n} \; .
\ee
\label{lem:scope}    
\end{lem}
provided $\gamma^{-k_n}>\delta$.

\noindent{\bf Proof:}  
If a curve traverses 
a cylinder of length $L$, width $(9/(\sqrt{\gamma})\, L$, then it also 
traverses a cylinder of width $(10/\sqrt{\gamma})\,L$ and length
$L/2$ centered at a line segment joining
discretized points in $L^\prime\, \Z^d$, 
provided $L^\prime\ \le L/\gamma$. 
The number of possible positions of such a cylinder
in a set of diameter $\ell$ is 
bounded above by $(\ell/L^\prime)^{2d}$. The number
of positions of  $n$ nested cylinders at scales $L_{k_1},
\dots , L_{k_n}$ is thus bounded by
\be
K_\Lambda\,\gamma^{2d\,k_1}\ \gamma^{2d\, (k_2-k_1)}\ldots 
\gamma^{2d\, (k_n-k_{n-1})}  \ \le\ K_\Lambda \gamma^{2d\, k_n}
\ee

Fix now a sequence $A_i$, $i=1,\dots , n$
of nested cylinders of length $L_{k_i}/2$ and width 
$(10/\sqrt{\gamma})\,L_{k_i}$. Let $\sigma$ be aspect ratio for
which ${\bf H2}$ holds with some $\rho<1$.  Cut each of the cylinders
into $\sqrt{\gamma}/(10\sigma)$
shorter cylinders of aspect ratio $\sigma$, and pick a maximal number
of well separated cylinders from this collection. Since
$A_{i+1}$ intersects at most two of the shorter  cylinders
obtained by subdividing $A_i$, the number of cylinders
in a maximal collection
is at least $n\left(\sqrt{\gamma}/(20\sigma)-2\right)$.
The probability of a curve traversing all of the $A_i$ is 
bounded above by the probability of crossing the shorter cylinders. 
Applying {\bf H2} gives
\be
\P_\delta\left(\begin{array}{c}
 A_1^\prime,\dots ,A_n^\prime \ \mbox{are crossed}\\
\mbox{by a curve in $\F_\delta(\omega)$} 
\end{array} \right)
\ \le\ K_1 e^{(K_2-K_3\sqrt{\gamma})\, n}
\ee 
Summing over the possible positions  and adjusting the constants
completes the proof.

\bigskip\noindent{\bf Proof of Theorem~\ref{thm:roughness}} \  We first
show that for each system of random curves in a compact set 
$\Lambda \subset \R^d$
satisfying the hypothesis {\bf H2}, 
there exist $m<\infty$ and $q<1$ such that for every $\gamma >  m$ 
\be 
Prob_{\delta}\left(  
\begin{array}{c}
	\mbox{ straight runs are $(\gamma,k_o)$-sparse} \\ 
	\mbox{ in $\Lambda$, down to scale $\delta$ } \\
\end{array} \right) \ \ge \ 1 - \frac{ q^{k_o} } {1- q} \ ;
\label{eq:probsparse} 
\ee 
in other words, the random variable given by
\be
k_{o;\delta}(\omega)\ =\ \inf\{ k\ge 0\ \mid\ 
\mbox{straight runs are $(\gamma,k_o)$-sparse down to scale 
$\delta$} \}
\label{eq:ko}
\ee
is stochastically  bounded as $\delta\to 0$.
To see this, note that for specified  $k$, 
\bea 
Prob\left(  
\begin{array}{c}
	\mbox{ there exist a nested sequence }  \\
	\mbox{ of straight runs on scales } \\
	\mbox{ $ k_1< \ldots < k_n = k$  with  $ k \ge n \ge k/2   $} 
\end{array} \right) \   & \le &  \  \sum_{ n=k/2}^{k} 
{k \choose n} \, K_1\,m^{2d\, k}\, e^{-K_2\sqrt{m}\, n} 
\nonumber \\ 
\ & \le & \  K_1\, (2m)^{2d\, k}\,e^{-(K_2\sqrt{m})\, k/2} \; .  
\eea 
Choosing $m$ large enough so that
\be
 q\ \equiv\ \, m^{2d}\, e^{-(K_2\sqrt{m}/2)}\ < \ 1 \; 
\ee 
and summing the geometric series over $k$ we obtain 
\eq{eq:probsparse}. 

As in the proof  of Theorem~\ref{thm:lower-det},
it follows with Lemmas~\ref{lem:fractal} and ~\ref{lem:cap-bound}
from \eq{eq:ko}
that all curves in a given configuration satisfy
the bound
\be
\CAP_{s;\delta}(\C) \ \ge \
 (\eps\,diam(\C))^{s} \left[ \gamma^{s\,k_{o;\delta}(\omega)}  \ + \
    \frac{ \beta}{1-\gamma^{s}/\beta}\right]^{-1} \; ,
\ee
with $m$ and $\gamma$ as above, $\beta=\sqrt{m(m+1)}$, and $s$
small enough so that $\gamma^s<\beta$. Choosing $\gamma$ sufficiently
close to $m$ we may take $s>1$, which proves the claim.
\hfill\blackbox

\startappendix 
 \vskip3cm 
\noindent {\large\bf APPENDIX} \vskip-0.6cm

\maappendix{Models with random curves} 
\label{sect:models}

In order to provide some context for the discussion of 
systems of random curves,
we present here a number of guiding examples.  
Familiarity with this material is not necessary for 
reading the work, however it does offer a better 
perspective both on the motivation and on the choice of 
criteria employed here.   
We start with some systems exhibiting 
the percolation transition.  

\masubappendix{Percolation models}

Among the simplest examples to present (for a review see 
\cite{SA,G}) 
is the independent bond percolation model on 
the cubic $d$ dimensional lattice, which we scale down to 
$\delta \Z^d$, $\delta << 1$. ``Bonds'' are pairs $b=\{x,y\}$ 
of  neighboring lattice sites.   Associated with  them  
are independent and identically distributed  random variables 
$n_b(\omega)$, with values in 
$\{0,1\}$.  The one-parameter family of probability 
measures is parametrized by:
\be
p \ = \ Prob\left( n_b  = 1 \right) \; .
\ee   
For a given realization, the bonds with $n_b(\omega)=1$ are referred 
to as occupied.  
The lattice decomposes into clusters of connected sites, with 
two sites regarded as connected if there is a path of occupied bonds 
linking them.   

For an intuitive grasp of the terminology one may think of the 
example in which the occupied bonds represent electrical conductors 
(of size $\delta << 1$) embedded randomly in an insulating medium.  
If a macroscopic piece 
of material with such characteristics is placed between two 
conducting plates which are maintained at different potentials, 
the resulting current
will be restricted to the macroscopic-scale clusters connecting the
two plates (the ``spanning clusters'').  

The model exhibits a phase transition.  Its simplest manifestation is
that the probability of there being an {\em infinite 
cluster} changes from $0$ for $p<p_c$, to $1$ for $p>p_c$.
The transition is also noticeable in finite volumes of macroscopic 
size: for $p<p_c$ the probability of observing a 
{\em spanning cluster} in $[0,1]^d$ is vanishingly small, 
whereas for  $p>p_c$ this probability is extremely close to $1$. 
In both cases the probabilities of the unlikely events 
decay as  $\exp{(-const.\ /\delta)}$, when $\delta \to 0$ at fixed
$p \ (\neq p_c)$.   

The generally believed picture in dimensions $2\le d<6$
is that for $p$ in the vicinity of the critical point 
($|p-p_c | = O(\delta^{1/\nu})$, macroscopic clusters do occur 
but are tenuous.  
Much of this is proven for 
$2D$ (\cite{R,SW,K-scaling}) though gaps in proof remain 
for $d>2$ (\cite{BCKS,Aiz-ISC}).
Typical configurations 
exhibit many {\em choke points},  where the change of the occupation 
status of a single bond will force a large scale
shift in the available connecting routes (\cite{R}), 
and possibly even break a connected cluster into two large 
components, as indicated in Figure~\ref{fig:chokepoint}.  
The clusters  are ``fractal'' in the sense that they exhibit 
fluctuating structure on many scales~\cite{Man}.  
This is the situation addressed in this work. 

For a given configuration of the model, we let 
${\cal F}_{\delta}(\omega)$  
stand for the collection of all the self-avoiding 
paths along the occupied bonds (possibly restricted to a 
specified subset $\Lambda \subset \R^d$).  
This random configuration of paths 
provides an explicit way of keeping track of the possible 
connecting routes within a given bond configuration.

One of the goals of this work was to establish that the 
description of the model in terms of a system of random curves 
(\cite{Aiz-IMA})  
remains meaningful even in the scaling limit ($\delta \to 0$).  
It may be noted that the alternate (and more common) 
description of the random configuration in terms of  the 
collection of {\em connected clusters}, is problematic 
in that limit.  
Clusters are naturally viewed as elements of the space 
of closed subsets of $\R^d$, with the distance provided 
by the Hausdorff metric.  
As long as $\delta \ne 0$ the two formulation of the model, 
as a system of random clusters or a system of random curves, 
are equivalent.  However 
the ubiquity of choke points renders the random cluster 
description insufficient for the scaling limit. 
(The  Hausdorff metric is not sensitive enough to pick up 
small differences, such as flips of individual bonds,
which may have a drastic effect on the available routes.)

It is expected that in the scaling limit the configurations 
of the connected 
paths in the critical bond percolation model are hard 
to distinguish from those arising from a number of 
other systems of different microscopic structure, e.g., 
percolation models where the conducting objects 
are randomly occupied sites of the lattice $\delta 
\Z^d$ (viewed as a subset of $\R^d$), or droplets 
of radius $\delta$ randomly 
distributed in $\R^d$.  The definition of 
${\cal F}_{\delta}(\omega)$ for the such models may require 
minor  adjustments, one in the notion of self-avoidance
and the other in the selection of the polygonal approximation.   
For the droplet model both are taken care of by 
restricting the attention to the polygonal paths joining 
centers of intersecting droplets 
{\em  which do not re-enter} any of the droplets.  
We form the set   
${\cal F}_{\delta}(\omega)$ as the collection of 
all such paths.  
 
In two dimensions, our hypotheses {\bf H1}  and {\bf H2} 
are satisfied by the independent bond, site, 
and droplet percolation models.  The $k=1$ 
case of the bounds (\ref{eq:suff-upper-1}) (with $\lambda(1) > 0 $)  
and (\ref{eq:suff-lower}), 
is a particular implications of the Russo-Seymour-Welsh  
theory  \cite{R,SW,Alex-RSW}. 
The statement that $\lambda(k) \to \infty$ follows by 
the van den Berg -- Kesten inequality \cite{vdBK}, 
which implies that for independent systems 
the probability of multiple crossings is dominated by the 
corresponding product of the probabilities of single events. 
(More detailed analysis implies that 
$\lambda(k)$ actually grows quadratically in $k$
\cite{Aiz-ISC,DupSal,Car97}.)  
The conditions {\bf H1} and {\bf H2} are  expected to  
hold also for other dimensions $d<6$, but not for 
$d>6$ (\cite{Aiz-ISC}).  
 
Thus, our general results imply the following 
statement, which was outlined in  ref.~\cite{Aiz-IMA}.

\begin{thm} 
In two dimensions, in each of the above mentioned percolation models, 
based on random bonds, random  sites, or random droplets, 
at the critical point all the non-repeating paths supported on the 
connected clusters within the compact region $[0,1]^{2}$
can be simultaneously parametrized by functions
$\gamma (t)$, $0\le t\le 1$, satisfying the H\"older 
continuity condition \eq{eq:holder-gen}.
The continuity constants $\kappa_{\eps;\delta}(\omega)$, which 
cover simultaneously all curves in $[0,1]^2$ 
remain stochastically bounded as $\delta \to 0$.  
(This holds for any $\eps >0$  as explained in 
Theorem~\ref{thm:regularity}). 

Furthermore, for each of these critical models 
the probability distribution of the random collection 
of curves  $\F_{\delta}(\omega)$ has a limit (in the sense of 
Theorem~\ref{thm:scaling}), at least for some sequence of 
$\delta_{n}\to 0$.  The limiting measure is supported on 
collections of curves whose Hausdorff dimensions satisfy 
\be
d_{min}\ \le  \ \mbox{dim\/}_{\cal H}(\C) \ \le d - \lambda(2) \; ,
\ee 
with some non-random $d_{min}>1$.
\end{thm} 
In fact, by similar reasoning we can also deduce 
the existence of a one-parameter 
family of such limits, corresponding to values of $p$ which deviate 
from $p_c$ by an amount which is scaled down to zero as 
$\delta \to 0$ 
 (in essense:  $p(\delta; t) = p_c + t \delta^{1/\nu}$).

The apparent universality of critical behavior 
leads one to expect that the scaling limits constructed here 
are common to the models listed above.
If so, then the limiting measures will have 
the full rotation and reflection symmetry of $\R^n$  
(and in two dimensions exhibit also {\em self-duality}).  
Remarkably, there is evidence for an 
even higher symmetry: conformal invariance 
(see \cite{LPS,Car92,Aiz-IMA,BS}), at the special point $p=p_c$, 
i.e., $t=0$ for the one-parameter family alluded to above.  
The mathematical derivation of such universality of 
the scaling limits, and of the conformal invariance at 
the critical point, form outstanding open problems.  

\masubappendix{Random spanning trees}

The regularity criteria presented here can also be verified
for a number of random spanning trees 
in two dimensions \cite{ABNW}.  Each is a translation-invariant
process describing a tree graph spanning a set of sites 
in $\R^d$ with neighboring sites spaced distances of order 
$\delta \ll 1$ apart.
\begin{itemize}
\item[MST]
  ({\em Minimal Spanning Tree}) \\
  The underlying graph is the 
  regular lattice $\delta \Z^d \subset \R^d$, with edges 
   connecting nearest neighbors. 
   Associated with the edges $b=\{x,y\}$ is a collection of 
   independent random {\em call numbers} (or {\em edge lengths}) 
   $u(b)$, with the uniform probability distribution in $[0,1]$.
   For a bounded region, $\Lambda \subset \R^d$, 
   the minimal spanning tree $\Gamma_{\delta;\Lambda}(\omega)$ is 
   the tree spanning the set $\Lambda \cap \delta \Z^d $  
   minimizing the {\em total edge length} 
   (i.e.\ the sum of the call numbers).  
\item[EST]
   ({\em Euclidean (Minimal) Spanning Tree})\\
    The vertices of the graph are generated as a random collection
    of points, with the Poisson distribution of density 
    $\delta^{-d}$ on $\Lambda$.
    $\Gamma_{\delta;\Lambda}$ is the covering tree graph 
    which minimizes the total (Euclidean) edge length.
\item[UST]
  ({\em Uniformly Random Spanning Tree}) \\
  The spanning tree $\Gamma_{\delta;\Lambda}$ is drawn uniformly 
   at random from the set of trees spanning the vertices in 
    $\Lambda\cap \delta\Z^d$ using the nearest neighbor edges.
\end{itemize}

In each of the above cases there is a well-defined limit 
\be 
\Gamma_{\delta}(\omega) \ = \ \lim_{\Lambda \nearrow {\R}^d} 
          \Gamma_{\delta,\Lambda}(\omega) 
\label{eq:tree}
\ee 
where  $\lambda$ is increased 
through a sequence which exhausts $\R^d$.  
(The restrictions 
of $\Gamma_{\delta,\Lambda}$ to compact subsets  
$\tilde \Lambda\subset\R^d$ are monotone decreasing in  
$\Lambda$  once 
$\Lambda\supset \tilde \Lambda$.) 
The limiting spanning tree is independent of the sequence 
of volumes and is translational invariant, in the stochastic 
sense.  

In general, the limit  $\Gamma_{\delta}(\omega)$ may be 
either a single tree or a collection of trees. 
For two dimensions it is known that each of MST, EST and UST 
almost surely consists of a single tree, with a single topological 
end, i.e. a single route to infinity 
(\cite{Pem91,Hagg95,BLPS98,Alex-MST,CCN,Alex-RSW}).  
The structure of UST changes from a tree to a forest
in dimensions $d>4$\cite{Pem91}, while MST and EST 
are expected to change similarly for $d>8$ \cite{New-Stein1, 
New-Stein2} (the transition may appear differently from 
the scaling limit perspective, \cite{ABNW}).

For each $n$-tuple of points 
$x_1, \ldots, x_n \in \R^d$, 
let $T^{(n)}_{x_1,\ldots,x_n}(\omega)$ be the tree subgraph of 
$\Gamma_{\delta}(\omega)$ with vertices corresponding to the closest 
$n$ sites in $\Gamma_\delta(\omega)$.
Our methods can be applied to the question, analogous to $Q1$ in the 
introduction: 
\begin{itemize}
{\em \item[Q.] Is there a limiting distribution for these graphs, as 
$\delta \to 0$?} 
\end{itemize}
To control the limit for $T^{(n)}_{x_1,\ldots,x_n}(\omega)$, 
one needs information on the curves supported on  
$\Gamma_{\delta}(\omega)$. 
This collection of curves forms the set  
${\cal F}^{(2)}_{\delta}(\omega)$ to which the analysis of this work 
may be applied.   

Random spanning trees provide striking examples of the phenomenon 
we encountered in critical percolation, that the formulation of the 
model in terms of random clusters is inadequate for the description 
of the scaling limit.  Here, the Hausdorff 
distance between any two different realizations (as subsets of 
$\R^d$) is~$\delta$, and hence the space of configurations seems to 
collapse to a single point.  That can be resolved by looking at 
the curves, as is done here.  
Let us add that the more complete description of the 
spenning trees requires the consideration of all the embedded finite 
trees, and that defines the object ${\cal F}_{\delta}(\omega)$ for 
those systems.  However, their study can be based on the analysis of 
the curves which provide the tree branches.

In contrast with independent percolation, 
spatially separated events are not independent for stochastic trees.  
Moreover, $\lambda(1) =0$, since any two vertices
are connected with probability one. 
Nevertheless,  the hypotheses {\bf H1} and {\bf H2} are valid 
(with $\lambda(2)>0$) for 
the three spanning tree processes
processes in $d=2$ dimensions \cite{ABNW}.  
Instrumental in the derivation are the relations of MST and 
EST with invasion percolation (studied in ref.~\cite{CCN}), and 
of UST with the
loop-erased random walk via the {\em Wilson algorithm} \cite{Wilson},
The latter relation permits to draw also non-trivial conclusions 
about the scaling limit of the loop-erased random walk (LERW) 
in d=2 dimensions. 

\masubappendix{ The frontier of Brownian motion.}

Yet another example of  a random curve is provided by 
the {\em frontier} of the two dimensional {\em Brownian motion} 
($\{\ b(t) \ | \ t \in [0,1], b(0)=0 \}$), 
abbreviated here as FBM. 
The frontier of a sample path is defined as 
the boundary of the unbounded connected component 
of the complement of the path in $\R^2$.

For the FBM, $\F(\omega)$ consists of a single curve.  
Its dimension has been considered in the literature: 
it is conjectured that $ dim\mbox{(FBM)}= 4/3$ (almost 
surely) \cite{Flory,Man}, the best rigorous bounds are 
$1.015\le dim\mbox{(FBM)} \le 1.475$ \cite{Law96,BL90}. 

Our  general results do apply to this 
example.  We shall not derive here the hypotheses 
{\bf H1}/{\bf H2}.  Let us note, however, that {\bf H2} 
is easy to establish by making use of the observation 
that the event depicted in Figure~\ref{fig:H2} requires that 
the Brownian path should hit each of the boxes  but 
not traverse it in the width direction. 
Thus, the mechanism behind 
our lower bound is similar in spirit to the 
earlier work of Bishop et.al.~\cite{BJPP}, 
in the reliance on the fact that Brownian paths 
move erratically.  The resulting upper bound, while not as 
tight an estimate of the dimension as that  of 
Burdzy and Lawler~\cite{BL90}, is expressed
as a bound on the tortuosity, and hence can be used to 
establish that FBM is parametrizable as a H\"older 
continuous curve. 


\masubappendix{The trail of three dimensional Brownian motion}

The trail of Brownian motion is the set of sites it visits 
in times $0\le t < \infty$.  In the transient case, $d>2$, 
the trail almost surely forms  a closed random set of Hausdorff 
dimension $2$.  Can it support curves of 
dimension arbitrarily close to $1$? 
In a recent work of Lawler \cite{Law3d} this question was 
answered negatively for the interesting case 
$d=3$, through analysis involving a number of results 
concerning the Brownian motion intersection exponent. 
Let us note that a negative answer can also 
be deduced from the general Theorem~\ref{thm:roughness}, 
since {\bf H2} is rather easy to establish (for $d>2$) 
within the  setup  relevant for this problem.

\masubappendix{Contour lines of random functions}

As the last example of a system of random lines let us 
mention contour lines of a random function.  Kondev and Henley 
\cite{KonHen} have considered the distribution of the level 
sets of a family of random functions defined on a lattice, 
$\phi_{\delta}(\omega): \delta  \Z^2 \rightarrow \R$.  
They present an interesting conjecture concerning the 
scale invariance for the distribution of  the loops 
bounding the connected  regions with $\phi(x)>\phi(0)$.
It would be of interest to see an extension of 
our analysis to such systems.

\noindent {\large \bf Acknowledgments\/} 
In the course of this work we have enjoyed stimulating 
discussions with numerous colleagues. 
In particular we would like to express our gratitude 
for useful comments to R. Langlands, S. Goldstein, 
C.M. Newman, and T. Spencer. 
We are especially grateful to  L. Berlyand for 
useful discussions of tortuosity lower bounds.  
This work was supported by the NSF Grant PHY-9512729 and, 
at the Institute for Advanced Study, by a grant from  
the  NEC Research Institute. 

\bigskip  \bigskip  
\addcontentsline{toc}{section}{References}


\begin{thebibliography}{10}

\bibitem{R}
L.~Russo, ``A note on percolation,'' {\em Zeit. Wahr.},  {\bf 43},  39 (1978).

\bibitem{SW}
P.~Seymour and D.~Welsh, ``Percolation probabilities on the square lattice,''
  in {\em {\em Advances in Graph Theory} Annals of Discrete Mathematics}
  (B.~Bollob{\'a}s, ed.), vol.~3, North Holland, 1978.

\bibitem{Alex-RSW}
K.~Alexander, ``The {RSW} theorem for continuum percolation and the {CLT} for
  {E}uclidean minimal spanning trees,'' {\em Ann. Appl. Probab.},  {\bf 6},
  466 (1996).

\bibitem{Aiz-ISC}
M.~Aizenman, ``On the number of incipient spanning clusters,'' {\em Nucl. Phys.
  B [FS]},  {\bf 485},  551 (1997).  
   \newblock {\tt http://xxx.lanl.gov/ps/cond-mat/9609240}.

\bibitem{ABNW}
M.~Aizenman, A.~Burchard, C.~Newman, and D.~Wilson, ``Scaling limits for
  minimal and random spanning trees in two dimensions,'' 1998 preprint.

\bibitem{BTW}
P.~Bak, C.~Tang, and K.~Wiesenfeld, ``Self--organized criticality: An
  explanation of $1/f$ noise,'' {\em Phys. Rev. Lett.},  {\bf 59},  381 
  (1987).

\bibitem{Bill}
P.~Billingsley, {\em Convergence of Probability Measures}.
\newblock John Wiley \& Sons, Inc., New York, 1968.

\bibitem{KZ93}
H.~Kesten and Y.~Zhang, ``The tortuosity of occupied crossings of critical
  percolation,'' {\em J. Statist. Phys.},  {\bf 70},  599 (1993).

\bibitem{Richardson}
L.~Richardson, ``The problem of contiguity,'' in {\em {G}eneral {Systems}
  {Y}earbook},  139, 1961.

\bibitem{Man}
B.~Mandelbrot, {\em The fractal geometry of nature}.
\newblock W. H. Freeman and Co., San Francisco, 1982.

\bibitem{DDC}
S.~Dunbar, R.~Douglass, and W.~Camp, ``The divider dimension of 
  the graph of a
  function,'' {\em J. Math. Anal. Appl.},  {\bf 167},  403 (1992).

\bibitem{Sta77}
H.~E. Stanley, ``Cluster shapes at percolation thresholds: 
  An effective cluster
  dimension and its connection with critical -- point exponents,'' {\em J.
  Phys. A},  {\bf 10},  L211 (1977).

\bibitem{Proh}
Y.~V. Prohorov, ``Convergence of random processes and limit theorems in
  probability theory,'' {\em Theor. Probability Appl.},  {\bf 1},  157,
  (1956).

\bibitem{Donsker}
M.~Donsker, ``An invariance principle for certain probability limit theorems,''
  {\em Mem. AMS},  {\bf 6} (1951).

\bibitem{Falconer}
K. J. Falconer, ``The geometry of fractal sets,''
  \newblock Cambridge University Press, 1985.

\bibitem{Erdos-Gillis}
P.~Erd{\H{o}}s and J.~Gillis, ``Note on the transfinite diameter,'' {\em
  Journal of the London Mathematical Society},  {\bf 12},  185 (1937).

\bibitem{SA}
D.~Stauffer and A.~Aharony, {\em Introduction to Percolation Theory}.
\newblock Taylor and Francis, London, 1994.

\bibitem{G}
G.~Grimmett, {\em Percolation}.
\newblock Springer-Verlag, 1989.

\bibitem{K-scaling}
H.~Kesten, ``Scaling relations for percolation,'' {\em Commun. Math. Phys.},
  {\bf 109},  109 (1987).

\bibitem{BCKS}
J.~Chayes, C.~Borgs, H.~Kesten, and J.~Spencer, ``Birth of the infinite
  cluster: finite size scaling in percolation.'' (In preparation).

\bibitem{Aiz-IMA}
M.~Aizenman, ``Scaling limit for the incipient spanning clusters,'' in {\em
  Mathematics of Materials: Percolation and Composites} (K.~M. Golden, G.~R.
  Grimmett, R.~D. James, G.~W. Milton, and P.~N. Sen, eds.), The IMA 
  Volumes in
  Mathematics and its Applications, Springer-Verlag 1998.
\newblock {\tt http://xxx.lanl.gov/ps/cond-mat/9611040}.

\bibitem{vdBK}
J.~van~den Berg and H.~Kesten, ``Inequalities with applications to percolation
  and reliability,'' {\em J. Appl. Prob.},  {\bf 22},  556 (1985).

\bibitem{DupSal}
B.~Duplantier and Saleur {\em Phys. Rev. Lett.},  {\bf 38},  2325 (1987).

\bibitem{Car97}
J.~Cardy, ``The number of incipient spanning clusters in two-dimensional
  percolation.'' 1997 preprint.

\bibitem{LPS}
R.~Langlands, P.~Pouliot, and Y.~Saint-Aubin, ``Conformal invariance in
  two-dimensional percolation,'' {\em Bull. AMS},  {\bf 30},  1 (1994).

\bibitem{Car92}
J.~Cardy, ``Critical percolation in finite geometries,'' {\em J. Phys. A},
  {\bf 25},  L201 (1992).

\bibitem{BS}
I.~Benjamini and O.~Schramm, ``Conformal invariance and {V}oronoi
  percolation.''
\newblock 1996 preprint.

\bibitem{CCN}
J.~Chayes, L.~Chayes, and C.~M. Newman, ``The stochastic geometry of invasion
  percolation,'' {\em Commun. Math. Phys.},  {\bf 101},  383 (1985).

\bibitem{Pem91}
R.~Pemantle, ``Choosing a spanning tree for the integer lattice uniformly,''
  {\em Ann. Probab.},  {\bf 19}, no.~4,  1559--1574 (1991).

\bibitem{Hagg95}
O.~{H\"aggstr\"om}, ``Random-cluster measures and uniform spanning forests,''
  {\em Stoch. Proc. Appl.},  {\bf 59},  267--275 (1995).

\bibitem{BLPS98}
I.~Benjamini, R.~Lyons, Y.~Peres, and O.~Schramm, ``Uniform spanning forests,''
  1998 preprint.
\newblock {\tt http://php.indiana.edu/\char126rdlyons/gz/usf.ps.gz}.

\bibitem{Alex-MST}
K.~Alexander, ``Percolation and minimal spanning forests in infinite graphs,''
  {\em Ann. Probab.},  {\bf 23},  87--104 (1995).


\bibitem{New-Stein1}
C.~M. Newman and D.~L. Stein, ``Spin-glass model with dimension-dependent
  ground state multiplicity,'' {\em Phy. Rev. Lett.},  {\bf 72},  2286--2289 
  (1994).

\bibitem{New-Stein2}
C.~M. Newman and D.~L. Stein, ``Ground state structure in a highly disordered
  spin glass model,'' {\em J. Stat. Phys.},  {\bf 82},  1113--1132 (1996).

\bibitem{Wilson}
D.~Wilson, ``Generating random spanning trees more quickly than the cover
  time,'' in {\em Proceedings of the {T}wenty-eighth {A}nnual {ACM} {S}ymposium
  on the {T}heory of {C}omputing ({P}hiladelphia, PA, 1996)},  296, {ACM},
  {N}ew {Y}ork, 1996.

\bibitem{Flory}
P.~Flory, ``The configuration of real polymer chains,'' {\em J. Chem. Phys.},
  {\bf 17},  303 (1949).

\bibitem{Law96}
G.~Lawler, ``The dimension of the frontier of planar {B}rownian motion,'' {\em
  Elect. Comm. in Probab.},  {\bf 1},  29 (1996).

\bibitem{BL90}
K.~Burdzy and G.~F. Lawler, ``Non-intersection exponents for {Brownian} paths,
  {Part II}: Estimates and applications to a random fractal,'' {\em Ann.
  Probab.},  {\bf 18},  981 (1990).

\bibitem{BJPP}
C.~J. Bishop, P.~W. Jones, R.~Pemantle, and Y.~Peres, ``The dimension of the
  {B}rownian frontier is greater than 1,'' {\em J. Funct. Anal.},  {\bf 143},
  309 (1997).

\bibitem{Law3d}
G.~Lawler, ``Hausdorff dimension of cut points for {B}rownian motion,'' {\em
  Elect. J. of Probab.},  {\bf 1},  no.2 (1996).

\bibitem{KonHen}
J.~Kondev and C.~Henley, ``Geometric exponents of contour loops on random
  Gaussian surfaces,'' {\em Phys. Rev. Lett.},  {\bf 74},  4580 (1995).

\end{thebibliography}

\end{document}